\documentclass{siamltex}

\usepackage{amsmath}
\usepackage{amssymb}
\usepackage{enumerate}

\usepackage{graphicx}
\usepackage{caption}
\usepackage{epstopdf}

\newcommand{\R}{\mathbb R}
\newcommand{\C}{\mathbb C}

\newcommand{\be}{\mathbf e}

\newcommand{\bg}{\mathbf g}

\newcommand{\bk}{\mathbf k}
\newcommand{\bm}{\mathbf m}
\newcommand{\bp}{\mathbf p}
\newcommand{\bq}{\mathbf q}
\newcommand{\br}{\mathbf r}
\newcommand{\bs}{\mathbf s}
\newcommand{\bt}{\mathbf t}
\newcommand{\bu}{\mathbf u}
\newcommand{\bv}{\mathbf v}

\newcommand{\bx}{\mathbf x}
\newcommand{\by}{\mathbf y}
\newcommand{\bz}{\mathbf z}

\newcommand{\bxk}{{\mathbf x}^{(k)}}

\newcommand{\byk}{{\mathbf y}^{(k)}}
\newcommand{\sk}{\sigma^{(k)}}
\newcommand{\vectornorm}[1]{\left\|#1\right\|}
\newcommand{\fA}{f(A)}
\newcommand{\BF}{\mathcal{F}}
\newcommand{\BW}{\mathcal{W}}
\newcommand{\K}{\mathcal{B}}
\newcommand{\Kp}{\widetilde{\K}}
\newcommand{\G}{\mathcal{G}}

\author{Sarah~W.~Gaaf\thanks{%
Department of Mathematics and Computer Science, TU Eindhoven, PO Box 513, 5600 MB, The Netherlands,
({\tt s.w.gaaf@tue.nl}). This author is supported by a Vidi research grant from the Netherlands Organisation
for Scientific Research (NWO).}
\and
Valeria~Simoncini\thanks{%
Dipartimento di Matematica, Universit{\`a} di Bologna,
Piazza di Porta S. Donato, 5, I-40127 Bologna, Italy ({\tt valeria.simoncini@unibo.it}).}
}
\title{Approximating the leading singular triplets of a large matrix function%
\thanks{Version of May 13, 2015.}}

\begin{document}
\bibliographystyle{siam}

\maketitle

\begin{abstract}
Given a large square matrix $A$ and a sufficiently regular function $f$ so that $f(A)$ is
well defined,
we are interested in the approximation of the leading singular values and corresponding
singular vectors of $f(A)$, and in particular of $\|f(A)\|$, where $\|\cdot \|$ is the matrix norm
induced by the Euclidean vector norm. Since neither $f(A)$ nor $f(A)v$ can
be computed exactly, we introduce and analyze
an {\it inexact} Golub-Kahan-Lanczos bidiagonalization procedure, where
the inexactness is related to the inaccuracy of the operations $f(A)v$, $f(A)^*v$.
Particular outer and
inner stopping criteria are devised so as to cope with the lack of a true residual.
Numerical experiments with the new algorithm on typical application problems are reported.
%
\end{abstract}

\section{Introduction}
Given a large $n\times n$ complex matrix $A$ and a sufficiently regular function $f$ so that
$f(A)$ is well defined,
we are interested in approximating the largest singular values and corresponding
singular vectors of the matrix function $f(A)$. This computation will also
give  an approximation to its 2-norm, namely,
$\|f(A)\|$, where $\|\cdot\|$ is the matrix norm induced by the Euclidean vector norm,
and it is defined as
\begin{eqnarray}\label{eqn:2norm}
\|f(A)\| = \max_{0\ne \bx\in{\mathbb C}^n} \frac{\|f(A) \bx\|}{\|\bx\|} .
\end{eqnarray}
In our presentation we will chiefly discuss this norm approximation because of
its interest in applications. However, we shall keep in
mind that the considered procedure allows us to also
determine both associated left and right singular vectors, and that
a group of singular triplets can be determined simultaneously.

 The problem of approximating the norm of a matrix
function arises in
the solution of stiff linear initial value problems \cite{Gear.81},\cite{Schmitt.83},
in the evaluation of derivatives and perturbations of matrix
functions, which arise for instance in electronic structure theory
\cite{Gil.2010},\cite{Rubensson.12},\cite{Jarlebring.Rubensson.12}, 
and in monitoring the magnitude of the inverse of distance matrices \cite{Baxter.94}.
In numerical linear algebra  the norm of matrix
polynomials may be used in the analysis of iterative procedures, and the
norm of rational matrix functions, and in particular of the transfer function may
give information on the sensitivity of the matrix itself to perturbations;
see, e.g.,  \cite{Trefethen.Embree.05},\cite{Choi.PhD.2013} and their
references.

If $A$ were normal, then the approximation could be stated in terms of an eigenvalue problem
in $A$.
Indeed, if $A=Q\Lambda Q^*$ is the eigendecomposition of $A$ with $Q$ unitary
and $\Lambda$ diagonal, then $f(A) = Q f(\Lambda) Q^*$ \cite{Horn.Johnson.91}, so
that the leading singular
values of $f(A)$ could be determined by a procedure that approximates the eigenvalues
of $A$.

The problem is significantly more challenging if $A$ is large and non-normal, since
there is no relation between eigenvalues and singular values that can be readily
exploited during the computation. Moreover, although  $A$ may be sparse, in general
$f(A)$ will be dense, and it cannot be computed explicitly. We are thus left
with procedures that use $f$ and $A$ by means of the action of
$f(A)$ to a vector $\bv$. The Lanczos bidiagonalization is among the most used
strategies for approximating selected singular triplets of a given matrix. Given a
matrix ${F}$, this procedure
generates a sequence of orthonormal vectors $\{\bv_1, \bv_2, \ldots \}$ and
$\{\bu_1, \bu_2, \ldots \}$ by alternating products of ${F}\bv$ and ${F}^*\bu$. In
our case, ${F} = f(A)$ and therefore these matrix vector products can be
approximately computed. In fact, since this computation is expensive,
 we shall consider an {\it inexact}
implementation of the Lanczos bidiagonalization process, where at each iteration
the action of $f(A)\bv$ and $f(A)^*\bu$ is approximated with some loose tolerance
by means of a projection
method. The problem of approximating $f(A)\bv$ has seen a great interest growth in the
past fifteen years, due to the emerging occurrence of this computation in many scientific
and engineering applications;
see, e.g., \cite{BenziSIREV.13},\cite{druskin:3760},\cite{FrommerMarch2006},%
\cite{Guettel.survey.13},%
\cite{Higham2008},\cite{Hochbruck.Ostermann.10},\cite{Hochbruck1999},
 and their references. For our purposes we shall use Krylov subspace methods for
approximating $f(A)^*\bu$ and $f(A)\bv$ at each iteration, equipped
with a cheap stopping criterion that may also be adapted to the outer current accuracy.
We shall show that the inexactness in the Lanczos bidiagonalization causes
the loss of the known symmetry structure of the process.
Nonetheless, as is the case in finite precision analysis \cite{Larsen.tr98},
orthogonality of the basis can be preserved, so that the recurrence
maintains its effectiveness.

If a rough approximation to $\|f(A)\|$ is the only quantity of interest, instead
of a group of singular triplets, then other approaches could
be considered. For instance, $f$ could be approximated by some other
more convenient functions, and then the resulting matrix function norm could
be more easily estimated. As an alternative, equivalent definitions of $f(A)$
could be used, from which the norm could also be estimated; or, the relation
of $\|f(A)\|$ with other norms or with some other spectral tool could be used; some
of these approaches are briefly recalled in section \ref{sec:known_methods}.
Methods in the mentioned classes, however, usually at most provide
the order of magnitude of the actual norm and are thus inappropriate
if more correct digits are needed.

This paper is organized as follows. Section~\ref{sec:known_methods} reviews some
methods available for the approximation of $\|f(A)\|$.
In section~\ref{sec: Lanczos bidiag} the
standard Lanczos bidiagonalization is recalled and the general notation used in this
paper is introduced. Section~\ref{sec: Inexact Lanczos bidiag} presents the
inexact Lanczos bidiagonalization procedure, including the details on the stopping
criteria in section \ref{sec:comput_stop}.
Section~\ref{sec:inner} discusses the approximation of the matrix function multiplication,
and a stopping criterion for its accuracy, while in section~\ref{sec:flex} a
stopping criterion in case an inner flexible strategy is analyzed.
In section~\ref{sec: spectral properties} we show how specific spectral properties allow us
to make a variable accuracy for the inner iteration feasible, which is finalized in
section~\ref{sec:variable}.
Section~\ref{sec: Practical implementation} focuses on the practical implementation and
the numerical results are presented in section~\ref{sec:expes}. We will conclude
with some discussion in section~\ref{sec: discussion}.

The following notation will be used throughout.
The vector $\be_i$ indicates the $i$th column of the identity matrix of a given dimension.
The conjugate transpose of a matrix $A$ will be denoted by $A^*$.
We will use the Matlab-like notation $[\bx; \by]$ to denote the column vector
$$
\begin{bmatrix} \bx \\ \by \end{bmatrix}, \quad \bx \in \C^{n_x}, \by\in \C^{n_y} .
$$
The Euclidean vector norm for vectors will be used, namely $\|\bx\| = (\sum_{i=1}^n |x_i|^2)^{\frac 1 2}$,
for $\bx\in\C^n$.
Unless explicitly stated, the induced matrix norm (\ref{eqn:2norm}) will be used for matrices.
For $A\in\C^{n\times n}$, {\rm spec}($A$) denotes the set of its eigenvalues, and
$W(A) = \{ z \in \C : z = (\bx^* A \bx)/(\bx^*\bx), \bx\in\C^n\}$ is its field of values.

\section{Available techniques for estimating the norm}\label{sec:known_methods}
While the Lanczos bidiagonalization is widely recognized as the method
of choice for approximating selected singular triplets of a large matrix,
if one is only interested in estimates
of $\|f(A)\|_2$ with $A$ non-Hermitian, then rather different procedures could also be used.
A simple approach consists of roughly estimating
$\| \cdot \|_2$ by using some other matrix norm. For instance,
$$
\|f(A)\|_2 \le \sqrt{n} \|f(A)\|_p, \quad p=1,\infty,
$$
where $\|f(A)\|_p$ is once again an induced norm \cite[p. 365]{Horn.Johnson.13}.
This bound is usually pessimistic, and it is clearly unsatisfactory for $n$ large.
The fact that for $A$ large the entries of $f(A)$ are not all readily available
 provides an additional challenge.

In the following we describe a few approaches available in the literature that
are tailored to the matrix function case.
Some of them first determine an explicit upper bound for the
norm, which only depends on scalar quantities.
The core computation will then be to determine a good approximation
to the obtained upper bound. 
The quality of the final estimate of $\|f(A)\|_2$ will thus depend both on
the sharpness of the initial upper bound and on the accuracy of the computation.
For general non-normal matrices the initial bound is often not very sharp,
limiting the quality of the overall estimation.
Finally, a computation-oriented estimation is the power method, which directly approximates
$\|f(A)\|_2$ as the square root of the largest eigenvalue
of $f(A)^*f(A)$. A more detailed list follows.

\begin{enumerate}
\item
Let $r(A)$ be the numerical radius of $A$, that
is $r(A) = \max\{|z|: z\in W(A)\}$, where $W(A)$ is the field of values of $A$.
Since  (see, e.g., \cite[Theorem 1.3-1]{Gustafson1997})
$$
r(A) \le \|A\|_2 \le 2\, r(A),
$$
by applying the bounds to $f(A)$ instead of $A$,
it is possible  to estimate $\|f(A)\|_2$ by means of $r(f(A))$; see, e.g.,
\cite{Watson.96},\cite{Mengi.Overton.05} for numerical methods to compute the numerical
radius of a given matrix.
A related special case is given by the exponential function, for which
the bound
\begin{eqnarray}\label{eqn:exp}
\|\exp(A)\|\le \exp(\alpha)
\end{eqnarray}
holds, where $\alpha$ is the largest eigenvalue
of the Hermitian part of $A$, that is of $\frac 1 2 (A+A^*)$ \cite[section 10.1]{Higham2008}.

\item If it is possible to find $K > 0$ and $\Omega\subset {\mathbb C}$ such that
$$
\|f(A)\|_2 \le K \|f\|_{\Omega},
$$
then it is sufficient to estimate $\|f\|_{\Omega}$; here $\|f\|_{\Omega}$ is the
$L_\infty$-norm of $f$ on $\Omega$. This can be done for instance
when $f$ is a polynomial, \cite{Crouzeix2007}, for which $K$ is known to be
less than 11.08 and conjectured to be equal to 2,
 and $\Omega$ coincides with the field of values of $A$. We refer
to the Ph.D. thesis of D. Choi \cite{Choi.PhD.2013},
for a discussion on the use if this bound when
$A$ is normal, or when $A$ is a contraction; see also
\cite{Trefethen.Embree.05} for a detailed analysis of this bound when
using pseudospectral information.
The computational intensive task is given by the determination of $\Omega$. If
$\Omega$ coincides with $W(A)$, then the cost of accurately approximating
$\Omega$ may be superior to that of approximating the single quantity
$\|f(A)\|_2$.

\item
This approach is in the same spirit as the one above, but with
a higher computational component.
For $\varepsilon>0$,
let $\sigma_{\varepsilon}(A) = \{z \in {\rm spec}(A+E): \|E\| < \varepsilon\}$
and assume that $f$ is analytic in $\sigma_{\varepsilon}(A)$.
If $L_\varepsilon$ denotes the length of the boundary
$\partial \sigma_{\varepsilon}(A)=\{ z\in\C: \|(zI-A)^{-1}\| =
\varepsilon^{-1}\}$, then by using the Cauchy integral expression for
$f(A)$ we obtain (see, e.g., \cite{Trefethen.Embree.05})
$$
\|f(A)\| \le \frac{L_\varepsilon}{2\pi\varepsilon}\|f\|_{\partial \sigma_{\varepsilon}} .
$$
Although the involved quantities may be easier to compute than in the previous
case, the dependence on $\varepsilon>0$ remains not fully controllable.

\item Using the relation $\|f(A)\|_2^2 = \lambda_{\max}( f(A)^*f(A))$ a
run of a few iterations of the power method can give an estimate to $\lambda_{\max}( f(A)^*f(A))$;
 see, e.g., \cite[Algorithm 3.19]{Higham2008} for an algorithm
specifically designed for the largest singular triplet.

\end{enumerate}

The power method is probably the most appealing approach among the
ones listed above. If a rough approximation is required,
typically to determine the order of magnitude, then its low cost provides
a satisfactory answer. However, if more than one digit of accuracy
is required, then the process may become slow. As for with $A$, the
stability in the computation may be highly influenced by the squaring; we refer
to section \ref{sec:expes} for an example of this well known phenomenon.

\section{Lanczos bidiagonalization}\label{sec: Lanczos bidiag}
We start by  recalling the Golub-Kahan bidiagonalization process in our context,
in terms of the matrix function
$f(A)$; then we will discuss how to actually obtain $f(A)$ times a vector.
Let $\bu_0 = 0$ and $\beta_1 = 0$, then for $i = 1,\ldots,m$ the following
recurrence relations define the Lanczos algorithm
\begin{equation}\label{eq: Lanczosbidiag}
\begin{array}{rcl}
\beta_{2i}\bu_i        &=& \fA\bv_i - \beta_{2i-1}\bu_{i-1}\\
\beta_{2i+1}\bv_{i+1}  &=& \fA^*\bu_i - \beta_{2i}\bv_{i}.\\
\end{array}
\end{equation}
The coefficients $\beta_j$ are computed so that the corresponding vectors $\bv$'s and
$\bu$'s have unit norm.
By collecting the two sets of vectors as
$U_m = \left[\bu_1, \ldots,\bu_m\right]$, $V_m = \left[\bv_1, \ldots,\bv_m\right]$, we
observe that $U_m^*U_m = I$, $V_m^*V_m = I$, $V_{m}^*\bv_{m+1} = 0$. Moreover, the
two recurrences can be compactly written as
\begin{eqnarray}
\fA V_m    &=& U_mB_m \nonumber \\
\fA^*U_m  &=& V_mB_m^* + \beta_{2m+1} \bv_{m+1}\be_m^*, \label{eqn:GK}
\end{eqnarray}
where $B_m$ is the following bidiagonal matrix
$$
B_m =
\begin{bmatrix}
\beta_{2} & \beta_3 & & \\
& \beta_{4} & \beta_5 & \\ &&\ddots&\ddots\end{bmatrix} \in \R^{m\times m} .
$$
It can be shown that the columns of $V_m$ span the Krylov subspace $\mathcal{K}_m(\fA^*\!\fA,\bv_1)$ and the columns of $U_m$ span the Krylov subspace $\mathcal{K}_m(\fA\fA^*,\fA\bv_1)$.
Define
$$
\K_{2m} =  \begin{bmatrix}0 & B_m \\ B_m^* & 0 \end{bmatrix} , \qquad \text{and} \qquad
\BW_{2m} = \begin{bmatrix} U_m & 0 \\ 0 & V_m \end{bmatrix} ,
$$
and
$$
\BF = \begin{bmatrix} 0 & f(A) \\ f(A)^* & 0 \end{bmatrix} .
$$
Then the recursion (\ref{eqn:GK}) can be rewritten as a standard Lanczos process,
in the more compact matrix notation (\cite[p.~178--186]{Cul02}, \cite[p.~448--449, p.~495]{GLo96})
\begin{equation}\label{eq: LBD with B 2}
\BF \BW_{2m} = \BW_{2m}\K_{2m} + \beta_{2m+1}\begin{bmatrix}0\\\bv_{m+1}\end{bmatrix} \be_m^*, \quad
\be_m\in\R^{2m}.
\end{equation}
For both even and odd $j$, the eigenvalues of $\K_j$ occur in $\pm$ pairs, with the exception of an extra extraneous zero eigenvalue in the odd case.
%
Within this setting, is it thus possible to
approximate the singular values of $\fA$
by the positive eigenvalues of $\K_{2m}$, or, equivalently, by the singular values of $B_{m}$.
In particular, for the largest singular value it holds that
(see \cite[Corollary $3.1.3$, Lemma $3.3.1.$]{Horn.Johnson.91}):
$$
\sigma_1(B_{j-1}) \le \sigma_1(B_{j}) \le \sigma_1(\fA), \qquad 2\le j\le m.
$$
There are several advantages of the Golub-Kahan bidiagonalization over the
simpler power method applied to $f(A)^*f(A)$, which are mainly related to the
fact that the eigenvalue squaring in this latter problem may lead to
severe  loss of information in the case very small or very large singular
values arise. In the inexact case the bidiagonal formulation also
allows us to better trace the inexactness during the whole
approximation process; this is discussed in the next section.

\section{Inexact Lanczos bidiagonalization}\label{sec: Inexact Lanczos bidiag}
When neither the explicit computation of the matrix $f(A)$ nor the
accurate operation $f(A)\bv$ (or $f(A)^*\bv$) are
feasible, then approximate computations must be performed,
resulting in an {\it inexact} Lanczos bidiagonalization procedure.
As a consequence, the recurrence (\ref{eqn:GK}) needs to be
significantly revised so as to aknowledge for the quantities that
are actually computed.

For a given $\bv$, the exact matrix-vector multiplication $f(A) \bv$ has to be replaced by
an inner procedure that approximates the resulting vector up to a certain accuracy.
The same holds for the operation $f(A)^*\bu$ for a given vector $\bu$.
For the sake of the analysis, at each iteration $k$ we shall formalize this
difference by writing, for some matrices ${C}_i$ and ${D}_i$,
\begin{eqnarray*}
\beta_{2i}\bu_i        &=& (\fA\bv_i+{C}_i\bv_i) - \beta_{2i-1}\bu_{i-1}\\
\beta_{2i+1}\bv_{i+1}  &=& (\fA^*\bu_i+ {D}_i \bu_i) - \beta_{2i}\bv_{i} ,
\end{eqnarray*}
where ${C}_i, {D}_i$ implicitly
 represent the perturbation induced by the approximate computations.
%
%
Since in general $\fA^* + {D}_i$ is no longer the conjugate transpose
of $\fA + {C}_i$, orthogonality of
a new vector $\bv_{m+1}$  has to be enforced by explicit
orthogonalization with respect to all previous vectors $\bv_i$, $1\le i \le m$.
The same holds for the vectors $\bu_i$, $i=1, \ldots, m$. Therefore, instead of one
bidiagonal matrix $B_m$ in the exact relation, we now obtain an
upper triangular matrix $M_m$ and an upper Hessenberg matrix $T_m$.
This leads to the following relations for the inexact (perturbed) Lanczos bidiagonalization:
\begin{eqnarray*}
(\fA + \mathfrak{C}_m) V_m   &=& U_mM_m\\
(\fA^* + \mathfrak{D}_m) U_m &=& V_m T_m + t_{m+1,m}\bv_{m+1}\be_m^*,
\end{eqnarray*}
where $\mathfrak{C}_m = \sum_{k=1}^{m}{C}_k\bv_k\bv_k^*$ and $\mathfrak{D}_m = \sum_{k=1}^m {D}_k\bu_k\bu_k^*$.
The matrices $V_m$ and $U_m$ are different from the matrices in the exact relation,
but they still have orthonormal columns.

The inexact Lanczos bidiagonalization can also be described using the notation of (\ref{eq: LBD with B 2}).
Define
\[\Kp_{2m} = \left[ \begin{array}{cc} 0 & M_m \\ T_m & 0 \end{array} \right], \qquad
\BW_{2m} = \left[ \begin{array}{cc} U_m & 0 \\ 0 & V_m \end{array} \right],\]
and the perturbation matrix 
$$
\G_{2m} = \begin{bmatrix} 0 & \mathfrak{C}_m \\ \mathfrak{D}_m & 0\end{bmatrix} \BW_{2m} =: {\cal E}_m \BW_{2m}.
$$
The perturbed relation thus becomes
\begin{equation}\label{eq: Inexact LBD K}
\BF \BW_{2m} + \G_{2m} = \BW_{2m}\Kp_{2m} +
t_{m+1,m}\left[\begin{array}{c} 0 \\ \bv_{m+1} \end{array}\right]\be_{m}^*, \quad \be_m \in\mathbb{R}^{2m},
\end{equation}
where
$$
\BF \BW_{2m} + \G_{2m}  = (\BF + {\cal E}_m) \BW_{2m}
=\begin{bmatrix}0 & \fA+\mathfrak{C}_m  \\ \fA^*+\mathfrak{D}_m & 0 \end{bmatrix} \BW_{2m}
 =: \widetilde\BF_{2m} \BW_{2m} .
$$
In contrast to the exact case, the space spanned by the columns of $\BW_{2m}$ is not a Krylov subspace.
However, when $t_{m+1,m}$ is small, this new space is close to an invariant subspace of the
perturbed matrix $\widetilde{\BF}_{2m}$, because then
$\widetilde{\BF}_{2m}\BW_{2m}\approx\BW_{2m}\Kp_{2m}$.
Notice the similarity of (\ref{eq: Inexact LBD K}) with equation (3.1) in \cite{Sim05}, which
shows that with this formulation, the inexact projection problem amounts to solving a
structured eigenvalue problem, where the original Hermitian matrix $\BF$ has been perturbed
by a structured non-Hermitian perturbation ${\cal E}_m$. The theory in \cite{Sim05} can then
be used to analyze and monitor the inexact computations, although the general results in
\cite{Sim05} should be carefully adapted to the new problem structure.

If ${\mathcal E}_m$ is small in norm, the eigenvalues of the {\it non-Hermitian} matrix
$\widetilde\BF_{2m}$ are small perturbations of the eigenvalues of the {\it Hermitian}
matrix $\BF$. Indeed, the eigenvalues of the
perturbed matrix $\widetilde{\BF}_{2m}$ lie in discs with radius
$\|{\cal E}_m\|$ and center the (real) eigenvalues of
$\BF$ (see, e.g., \cite[section IV, Theorem 5.1]{StS90}). Therefore, for small
perturbations in the computations, the eigenvalues of the symmetric matrix $\BF$ will be
perturbed accordingly.
On the other hand, in the following we shall consider the case when $\|{\mathcal E}_m\|$ is
larger than usually allowed by a perturbation analysis argument, therefore different
strategies need to be devised to ensure good approximations to the wanted eigenvalues of $\BF$.

Following the standard procedure of the exact case, we should consider the matrix
$\Kp_{2m}$ to approximate the largest eigenpairs of $\widetilde{\BF}_{2m}$, and according
to the discussion above, of $\BF$. Due to the non-Hermitian structure of $\Kp_{2m}$, however,
there are different matrices that can provide the sought after singular value information, namely
the matrix $\Kp_{2m}$ itself, and the two distinct
matrices $T_m$ or $M_m$. The last two matrices yield approximations to the corresponding triplets
of $f(A)+\mathfrak{C}_m$ and $f(A)^* + \mathfrak{D}_m$.
%
The following bound between the largest eigenvalue of $\widetilde\BF$ and
the largest singular values of $T_m$ and $M_m$ shows that all these quantities
can be easily related. Let $\bq = [\bx; \by]$.
Using $\|\bx\|\|\by\| \le \frac{1}{2} \left(\|\bx\|^2+\|\by\|^2\right)$ we obtain\footnote{Another
bound can be obtained for the geometric mean, that is $|\theta| \le \sqrt{\sigma_1(M_m)\sigma_1(T_m)}$.}
\begin{eqnarray*}
|\theta| &\le& \displaystyle\max_{\bq\ne 0}\left|\frac{\bq^*\Kp_{2m}\bq}{\bq^*\bq}\right|
= \displaystyle\max_{[\bx;\by]\neq 0} \left|\frac{\bx^*M_m\by + \by^*T_m\bx}{\bx^*\bx + \by^*\by}\right|\\
                  &\le&  \displaystyle\max_{[\bx;\by]\neq 0}\frac{\|\bx\|\|M_m\by\|}{\|\bx\|^2 + \|\by\|^2} +
 \max_{[\bx;\by]\neq 0}\frac{\|\by\|\|T_m\bx\|}{\|\bx\|^2 + \|\by\|^2}\le \frac{1}{2}(\sigma_1(M_m) + \sigma_1(T_m)).
\end{eqnarray*}
If the inexactness of the bidiagonalization is very large,
$M_m$ and $T_m^*$ are very different from each other. In this case,
the leading singular values of these two matrices - and thus their mean - may be
significantly larger than the biggest (in modulo) eigenvalue of $\Kp_{2m}$, since they are
related to the numerical radius of $\Kp_{2m}$, rather than to its spectrum.
This motivated us to use the eigenvalues of $\Kp_{2m}$ in the approximation, rather
than the singular values of its blocks. Moreover, working with $\Kp_{2m}$ made the
analysis of the relaxed strategy particularly convenient, since known results on
relaxed eigenvalue computation could be exploited.

\subsection{A computable stopping criterion}\label{sec:comput_stop}
In this section we analyze a strategy for monitoring the convergence of the
inexact bidiagonal iteration. As it is common to other inexact processes, the true
problem residual is inaccessible as soon as inexactness takes place. However,
on the one hand
some stopping criterion needs to be introduced to exit the process.
On the other hand, it is unclear whether the computed approximations are still
meaningful for the original problem, since they were computed with significantly
modified data.

Let $(\theta,\bq)$ be an eigenpair of $\Kp_{2m}$, where $\bq$ is a unit vector.
As the iterations proceed,
$(\theta, \BW_{2m}\bq)$ tends to approximate an eigenpair of $\widetilde\BF$.
We would like to ensure that $(\theta, \BW_{2m}\bq)$  also tends to an
eigenpair of $\BF$.
To monitor the convergence of $\theta$ and to define a stopping criterion
for the outer iteration, the residual is used.
We call $\BF \BW_{2m}\bq - \theta \BW_{2m}\bq$ the \emph{true residual}, which is not available,
since $\BF$ cannot be applied exactly.
We thus introduce the \emph{computed residual}, which is the residual of the
actually computed quantities,
namely (see (\ref{eq: Inexact LBD K}))
\[
\br_{2m} := \widetilde \BF \BW_{2m}\bq  - \theta \BW_{2m}\bq =
t_{m+1,m}\begin{bmatrix}0 \\ \bv_{m+1}\end{bmatrix}
\be_{m}^*\bq, \quad (\be_m \in\mathbb{R}^{2m}).
\]
In the sequel, we shall use the following obvious inequality to
estimate the true residual norm:
\begin{equation}\label{eq: true comp gap}
\|\BF \BW_{2m}\bq - \theta \BW_{2m}\bq \| \le \|\br_{2m}\| +
\|(\BF \BW_{2m}\bq - \theta \BW_{2m}\bq) - \br_{2m}\|,
\end{equation}
where
$\|(\BF \BW_{2m}\bq - \theta \BW_{2m}\bq) - \br_{2m}\|$ is the gap between the computed and
the true residuals, in short the ``residual gap''. If this gap can be imposed to be small, then
the computed residual will give an estimate for the true residual.
In this case, convergence can be monitored by only using the (available)
computed residual, and the following relative stopping criterion can be used:
\begin{equation}\label{eqn:outer_criterion}
{\rm if} \quad \frac{|t_{m+1,m}\be_m^*\bq|}{|\theta|} < \varepsilon_{out} \quad {\rm then} \quad {\rm stop}
\end{equation}
for some outer tolerance $\varepsilon_{out}$, where $\theta$ is the largest
(in modulo) eigenvalue.
Finally, as the computed residual norm goes to zero, the quantity
$\|(\BF \BW_{2m}\bq - \theta \BW_{2m}\bq) - \br_{2m}\|$ will tend to dominate
again, playing the role of the final attainable accuracy level.

%


To see how we can impose the residual gap to be small, and recalling
the definition of $\G_{2m}$,  we first consider a more convenient expression
for $\G_{2m}\bq$, with $\bq = [\bx;\by]$, that is
$$
\G_{2m}\bq = \begin{bmatrix} \mathfrak{C}_m V_m \by\\ \mathfrak{D}_m U_m  \bx \end{bmatrix}
=: \begin{bmatrix} G_m^{(1)} \by\\ G^{(2)}_m   \bx \end{bmatrix}.
$$
Let $G_m^{(i)} = [\bg_1^{(i)}, \ldots, \bg_m^{(i)}]$. Then
\begin{eqnarray}
&&\|(\BF \BW_{2m}\bq - \theta \BW_{2m}\bq) - \br_{2m} \|  = \|\G_{2m}\bq\|
=\left \| \begin{bmatrix} G_m^{(1)} \by\\ G^{(2)}_m   \bx \end{bmatrix}\right \|
=\left \| \sum_{j=1}^m \begin{bmatrix} \bg_j^{(1)} \be_j^* \by \\ \bg_j^{(2)} \be_j^* \bx\end{bmatrix}\right \|
\nonumber \\
&& \qquad \le\sum_{j=1}^m \left \| \begin{bmatrix} \bg_j^{(1)} \be_j^* \by \\ \bg_j^{(2)} \be_j^* \bx\end{bmatrix}\right \|
= \sum_{j=1}^m (\|\bg_j^{(1)}\|^2\, |\be_j^* \by|^2 +\| \bg_j^{(2)}\|^2\,
|\be_j^* \bx|^2)^{\frac 1 2}. \label{eqn:G}
\end{eqnarray}
The vectors $\bg_j^{(i)}$, $i=1,2$, implicitly carry the error
caused by the inexact computation of $\fA\bv$ and $\fA^*\bu$, respectively, in the inner iteration.
If every term of this sum is small, the computed residual will be close
to the true residual. The following Lemma states how the inaccuracy
in the matrix-vector products relates to the residual gap;
its proof is based on the corresponding result in \cite{Sim05}, however the structure
is exploited so as to have a dependence with respect to $m$ instead of $2m$, the size
of $\widetilde{\cal B}_{2m}$.

\begin{lemma}\label{lemma: res bound}
Assume that $m$ iterations of the inexact Lanczos bidiagonalization process have been taken.

If $\|\bg_j^{(1)}\|, \|\bg_j^{(2)}\| < \frac{1}{m}\varepsilon$ for $1\le j \le m$, then
$\|(\BF \BW_{2m}\bq - \theta \BW_{2m}\bq) - \br_{2m}\| < \varepsilon$.
\end{lemma}

{\it Proof.} From $\|\bq\|=1$ with $\bq=[\bx;\by]$
it follows that $\|[\be_j^* \bx; \be_j^* \by]\|\le 1$.
From (\ref{eqn:G}) we obtain
\begin{eqnarray*}
\|(\BF \BW_{2m}\bq - \theta \BW_{2m} \bq) - \br_{2m}\|
&\le &
\sum_{j=1}^m (\|\bg_j^{(1)}\|^2\, |\be_j^* \by|^2 +\| \bg_j^{(2)}\|^2\,
|\be_j^* \bx|^2)^{\frac 1 2} \\
& < & \sum_{j=1}^m \frac{1}{m}\varepsilon (|\be_j^* \by|^2+ |\be_j^* \bx|^2)^{\frac 1 2} \le
\frac{1}{m}\varepsilon \sum_{j=1}^m 1 = \varepsilon . \qquad \endproof
\end{eqnarray*}

This result shows that if $\varepsilon$ is sufficiently small, then
the residual gap will stay below the computed residual norm until
convergence. In our experiments, $m$ will play the role of the maximum
number of Lanczos bidiagonalization iterations, which is usually set
to a number between $50$ and $500$.


\section{Approximation of $\fA\bv$ and $\fA^*\bu$ and a computable inner stopping criterion}
\label{sec:inner}
The performance of the inexact Lanczos bidiagonalization process depends on the
approximation accuracy of the matrix-vector products $\fA\bv$ and $\fA^*\bu$.
Due to the size of $A$, we consider approximating
this quantity by means of a projection-type iterative method as follows; we limit our discussion to $\fA\bv$,
and a corresponding procedure can be used for $\fA^*\bu$.
Starting with the unit vector $\bv$ and the matrix $A$, we construct a sequence of approximation subspaces
${\cal K}_k$ of ${\mathbb R}^n$,
$k=1, 2, \ldots,$ and define the matrix $P_k = [\bp_1 , \bp_2, \bp_3, \ldots,\bp_k] \in\C^{n\times k}$, whose
orthonormal columns span the subspace, and $\bv = \bp_1=P_k \be_1$, in a way so that
the spaces are nested, that is ${\cal K}_k \subseteq {\cal K}_{k+1}$.
Typical such choices are Krylov and rational Krylov subspaces \cite{Higham2008},\cite{Guettel.survey.13}.
The desired approximation is then obtained as
\[
\fA\bv \approx P_kf(H_k)\be_1, \qquad H_k = P_k^* A P_k.
\]
For small $k$, the reduced non-Hermitian
matrix $H_k$ has small size, so that $f(H_k)$ can be computed efficiently by
decomposition-type methods \cite{Higham2008}.

Our stopping criterion of this approximation process is based on an estimation of the error norm, and
it uses an approach previously introduced in \cite[Proposition~2.2]{KnizhnermanSimoncini2010}.
\begin{proposition}{\rm \cite[Proposition 2.2]{KnizhnermanSimoncini2010}}
Assume that $k+j$ inner iterations have been executed.
Let $\bz_{k+j} = P_{k+j}f(H_{k+j})\be_1$ be an approximation to $\fA\bv$ and define
$\omega_{k+j} = \|\bz_{k+j} - \bz_k\| / \|\bz_k\|$.
If $\|\fA\bv - \bz_{k+j}\| \ll \|\fA\bv - \bz_{k}\|$ and $\|\fA\bv\| \approx \|\bz_{k}\|$, then
\begin{eqnarray}\label{eqn:error_estimate}
\|\fA\bv - \bz_k\| \approx \frac{\omega_{k+j}}{1 - \omega_{k+j}}\,\|\bz_k\|.
\end{eqnarray}
\end{proposition}

The result in (\ref{eqn:error_estimate})
shows that after $k+j$ iterations it is possible to provide an estimate of
the error norm at iteration $k$. Therefore, we introduce the following
stopping criterion for the approximation of $f(A) \bv$:
\[
{\rm if} \quad \frac{\omega_{k+j}}{1-\omega_{k+j}} \le \varepsilon_{in} \quad {\rm then} \quad {\rm stop}
\]
for some inner tolerance $\varepsilon_{in}$.
The accuracy of the inner iteration will influence the final accuracy of the inexact Lanczos bidiagonalization.
In the notation of the previous section, once the stopping criterion is satisfied,
we have thus derived the following estimate for the
perturbation occurring in the Lanczos step,
$$
\| f(A)\bv_k - \bu_k  \| = \| C_k\bv_k\| \approx \varepsilon_{in} \| \bu_k\| ;
$$
Note that here $\| C_k\bv_k\| = \|\bg_k^{(1)}\|$, with the notation in (\ref{eqn:G}).
An analogous relation holds with respect to $f(A)^*\bu_k$ and thus $\|\bg_k^{(2)}\|$.
We stress here that, since the approximation process changes at each iteration,
the quantity $\| C_k\bv_k\|$ will vary as the Lanczos bidiagonalization
proceeds.
The threshold itself
may vary during the Lanczos iteration, so that $\varepsilon_{in}= \varepsilon_{in}^{(k)}$.
As experienced with other eigenvalue and linear system problems, $\varepsilon_{in}^{(k)}$
may even be allowed to grow during the iteration, without significantly affecting
the overall process. This is discussed in the next section.

\section{Relaxing the inner solution accuracy}\label{sec:flex}
The bound in (\ref{eqn:G}) on the residual gap suggests that the accuracy
on the inner solution approximation can be relaxed as convergence takes place.
Indeed, following similar strategies in \cite{Simoncini2003c},\cite{Sim05},\cite{BoF00}, we observe that
it is the product $\|\bg_j^{(1)}\|\, |\be_j^*\by|$ in (\ref{eqn:G})
 that needs to be small, and not each factor, to ensure a small gap; the same for
$\|\bg_j^{(2)}\|\,|\be_j^*\bx|$.
Therefore, if $|\be_j^*\by|$ is sufficiently small, indicating that
the $(m+j)$th component of the eigenvector $\bq$ is small, $\|\bg_j^{(1)}\|$ is allowed to
be larger, and the required accuracy of $\varepsilon$ can still be achieved.
This induces a variable (possibly growing) accuracy in the inner iteration, which
drives the size of $\|\bg_j^{(1)}\|$.
In the following we shall first show that 
the quantities $|\be_j^*\by|$ and $|\be_j^*\bx|$  do tend to decrease as
the approximation improves. We then derive a computable expression for the
variable stopping tolerance in the approximation of $f(A)\bv$ and $f(A)^*\bu$ at each iteration of
the resulting ``relaxed'' Lanczos bidiagonalization process.
This strategy may be convenient
in case the cost approximating $f(A)\bv$ and $f(A)^*\bu$
is very high, as is the case for instance if an accurate approximation to the
leading singular triplets is requested.

\subsection{Spectral properties of the approximate singular triplets}\label{sec: spectral properties}
To ensure that the magnitude of $\|\bg_j^{(k)}\|$, $k=1,2$, can be relaxed in the bound (\ref{eqn:G}), we
need to verify that $|\be_j^*\bx|$ and $|\be_j^*\by|$
become small as convergence takes place. This fact
has been verified in the eigenvalue setting in \cite{Sim05}, however the peculiar structure
of the Lanczos bidiagonal recurrence requires the validation  of the results in \cite{Sim05}.
To this end,
we first define the submatrix of $\Kp_{2m}$ of size $2k$ as
\[
\Kp_{2k} = \begin{bmatrix} 0 & M_k \\ T_k & 0 \end{bmatrix} ,
\]
where $M_k$, $T_k$ are the leading portions of the corresponding $m\times m$ matrices.
Let $(\theta^{(2k)}, \bq^{(2k)})$ be an eigenpair of
$\Kp_{2k}$, where $\bq^{(2k)} = [\bx;\by]$ has unit norm,  and $\bx,\by\in\C^{k}$.
Further, let
\begin{equation}\label{eq: q tilde}
\widetilde{\bq} = \begin{bmatrix}\bx \\ 0\\ \by \\ 0 \end{bmatrix},
\end{equation}
where the $0$-vectors have length $m-k$, and define
$\mathcal{X} = [\widetilde{\bq}, Y]$, where $Y$ is chosen such that
$\mathcal{X}$ is unitary. Define $\underline{\Kp}_{2m} = Y^*\Kp_{2m}Y\in\C^{(2m-1)\times(2m-1)}$.
The following result shows that under certain hypotheses some of the components of the approximate
eigenvectors do tend to zero as convergence takes place.
Its proof is analogous to that of \cite[Prop.~2.2]{Sim05}, and can be found in
 the appendix.

\begin{proposition}\label{prop:tau}
Let $(\theta^{(2k)},\bq^{(2k)})$ be an eigenpair of $\Kp_{2k}$,
and $\widetilde{\bq}$ be as defined in \emph{(\ref{eq: q tilde})}.
 Let $\bs_{2m}^* = \widetilde{\bq}^*\Kp_{2m} - \theta^{(2k)}\widetilde{\bq}^*$, $\delta_{2m,2k} = \sigma_{\text{min}}(\underline{\Kp}_{2m} - \theta^{(2k)}I) > 0$,
and $\br_{2k} = t_{k+1,k}\left[\begin{array}{c}0\\ \bv_{k+1}\end{array}\right]\be_k^*\bq^{(2k)}$. If
\[
\vectornorm{\br_{2k}} < \frac{\delta_{2m,2k}^2}{4\vectornorm{\bs_{2m}}},
\]
then there exists a unit norm eigenvector $\bq = [\bx_1;\bx_2;\by_1;\by_2]$
of $\Kp_{2m}$ with $\bx_1,\by_1\in\C^k$, $\bx_2,\by_2\in\C^{m-k}$, such that
\[
\vectornorm{\begin{bmatrix}\bx_2 \\ \by_2 \end{bmatrix}} \le \frac{\tau}{\sqrt{1+\tau^2}},
\]
with $\tau\in\R$, $0\le \tau < 2\frac{\vectornorm{\br_{2k}}}{\delta_{2m,2k}}$. Moreover, if $\theta$ is the eigenvalue associated with $\bq$, we have
\begin{equation}\label{eqn:thetas}
|\theta - \theta^{(2k)}| \le \|\bs_{2m}\|\tau.
\end{equation}
\end{proposition}

This proposition states that if after $k\le m$ iterations of Lanczos bidiagonalization
the computed residual $\|\br_{2k}\|$ is sufficiently small, then there exists
 an eigenvector of $\Kp_{2m}$ such that some of its components
are bounded correspondingly. These are precisely the components that
allow us to relax the accuracy in the inner iteration.
Note that $\delta_{2m,2k}$ gives an indication of the distance between the spectrum of
$\underline{\Kp}_{2m}$ and $\theta^{(2k)}$. It should be kept in mind that for non-normal matrices,
the value of $\delta_{2m,2k}$ may be much smaller \cite[Example 2.4, p. 234]{StS90}.
On the other hand, since $\widetilde {\cal B}_{2m}$ is a perturbation to a Hermitian matrix,
the quantity $\bs_{2m}$ is an approximate residual for $(\theta^{(2k)}, \bq^{(2k)})$ as an eigenpair
of $\widetilde {\cal B}_{2m}$, and thus it will be small as $m$ grows.
As a consequence, the condition in the theorem is likely to be satisfied, and the
eigenvalue error (\ref{eqn:thetas}) may be much smaller than $\tau$.

\subsection{Variable accuracy in the inner approximation}\label{sec:variable}
In this section we show that
relaxation in the inner accuracy at step $k\le m$ is
possible if there exists an eigenpair $(\theta^{2(k-1)}, \bq^{2(k-1)})$ of $\Kp_{2(k-1)}$ such that
\begin{eqnarray}
\|\br_{2(k-1)}\| &<& \frac{\delta_{2m,2(k-1)}^2}{4\|\bs_{2m}\|}, \label{restriction 1} \\
\forall \theta_j \in \Lambda(\Kp_{2m}), & &
 \quad \theta_j \neq \theta, \qquad |\theta_j - \theta^{2(k-1)}| > 2\frac{\|\bs_{2m}\|\|\br_{2(k-1)}\|}{\delta_{2m,2(k-1)}} .  \label{restriction 2}
\end{eqnarray}
The first condition (\ref{restriction 1}) ensures that there exists an
eigenvector $\bq$ of $\widetilde{\cal B}_{2m}$ whose specified components
are small, according to Proposition \ref{prop:tau}. Let $\theta$ be the eigenvalue
associated with this $\bq$. The
second condition, (\ref{restriction 2}), guarantees
that the eigenvalue $\theta^{2(k-1)}$ of $\Kp_{2(k-1)}$ is a perturbation of the
eigenvalue $\theta$ of $\Kp_{2m}$, which is the final approximation to the original problem.
The following theorem states how the use of a variable accuracy will still guarantee
a small residual gap, and hence yields a true residual with an accuracy which
is bounded by the accuracy of the gap, in agreement with (\ref{eq: true comp gap}).

\begin{theorem}
Assume $m$ inexact Lanczos bidiagonalization iterations are carried out.
Let $(\theta,\bq)$ be an eigenpair of $\Kp_{2m}$, where $\theta$ is
simple and $\|\bq\| = 1$. Given $0<\varepsilon_{out} \in \mathbb{R}$, with the
notation of (\ref{eqn:G}) assume that for $k = 1,\ldots,m$
\begin{equation}
\|\bg_k^{(1)}\|,\|\bg_{k}^{(2)}\| \le \left\{
\begin{array}{ll}
\frac{\delta_{2m,2(k-1)}}{2m\|\br_{2(k-1)}\|}\varepsilon_{out} &
\begin{array}{l}
\text{if $k>1$, and there exists $(\bq^{2(k-1)},\theta^{2(k-1)})$} \\
\text{of $\Kp_{2(k-1)}$ satisfying \emph{(\ref{restriction 1})} and \emph{(\ref{restriction 2})},}
\end{array} \\
&\\
\frac{1}{m}\varepsilon_{out} & \begin{array}{l} \text{otherwise.}
 \end{array}
 \end{array}
 \right.
\end{equation}
Then $\|(\BF \BW_{2m}\bq - \theta \BW_{2m}\bq) - \br_{2m}\| \le \varepsilon_{out}$.
\end{theorem}

{\it Proof.}
This proof is analogous to the proof of
Theorem 3.1 in \cite{Sim05}. Suppose that at the $(k-1)$th iteration there exists an
eigenpair $(\theta^{2(k-1)}, \bq^{2(k-1)})$ of $\Kp_{2(k-1)}$ satisfying the
conditions $(\ref{restriction 1})$ and $(\ref{restriction 2})$.
This implies that $\theta^{2(k-1)}$ is a perturbation of the considered
eigenvalue $\theta$ of $\Kp_{2m}$, since $\theta$ is the only eigenvalue of $\Kp_{2m}$ such that
\[
|\theta - \theta^{2(k-1)}| \le  2\frac{\|\bs_{2m}\|\|\br_{2(k-1)}\|}{\delta_{2m,2(k-1)}}.
\]
Let $\mathcal{K}\subset \{1,\ldots,m\}$ be defined such that for each $k\in\mathcal{K}$ there exists a eigenpair $(\bq^{2(k-1)}, \theta^{2(k-1)})$ of $\Kp_{2(k-1)}$ satisfying the conditions $(\ref{restriction 1})$ and $(\ref{restriction 2})$.
Then, similar to the reasoning in the proof of Lemma \ref{lemma: res bound} and using (\ref{eqn:G}),
\begin{eqnarray*}
\|(\BF \BW_{2m}\bq &-& \theta \BW_{2m} \bq) - \br_{2m}\| =
 \|\G_{2m}\bq\| \le
\sum_{k=1}^m (\|\bg_k^{(1)}\|^2\, |\be_k^* \by|^2 +\| \bg_k^{(2)}\|^2\,
|\be_k^* \bx|^2)^{\frac 1 2} \\
  &\le&
\sum_{k\in\mathcal{K}}
(\|\bg_k^{(1)}\|^2\, |\be_k^* \by|^2 +\| \bg_k^{(2)}\|^2 |\be_k^* \bx|^2)^{\frac 1 2} +
\sum_{\substack{k\notin \mathcal{K},\\k\le m}}
(\|\bg_k^{(1)}\|^2\, |\be_k^* \by|^2 +\| \bg_k^{(2)}\|^2 |\be_k^* \bx|^2)^{\frac 1 2} \\
&\le&
\sum_{k\in\mathcal{K}}
 \frac{\delta_{2m,2(k-1)}\varepsilon_{out}}{2m\|\br_{2(k-1)}\|}
(|\be_k^* \by|^2 +|\be_k^* \bx|^2)^{\frac 1 2} +
\sum_{\substack{k\notin \mathcal{K},\\k\le m}}
\frac{\varepsilon_{out}}{m}
(|\be_k^* \by|^2 +|\be_k^* \bx|^2)^{\frac 1 2} \\
&\le&
\sum_{k\in\mathcal{K}}
 \frac{\delta_{2m,2(k-1)}\varepsilon_{out}}{2m\|\br_{2(k-1)}\|}
2 \frac{\|\br_{2(k-1)}\|}{\delta_{2m,2(k-1)}} +
\sum_{\substack{k\notin \mathcal{K},\\k\le m}}
\frac{\varepsilon_{out}}{m} \\
&= & \displaystyle\frac{|\mathcal{K}|}{m}\,\varepsilon_{out} +
\frac{m - |\mathcal{K}|}{m}\,\varepsilon_{out}  = \varepsilon_{out} \qquad
\endproof
\end{eqnarray*}

\begin{table}[htbp]\label{algorithm1}
{\footnotesize
\noindent \vrule height 0pt depth 0.5pt width \textwidth \\
{\bfseries Algorithm 1:} Inexact Lanczos bidiagonalization\\[-3mm]
\vrule height 0pt depth 0.3pt width \textwidth \\
{\bf Input: } $A\in{\mathbb C}^{n\times n}$ non-Hermitian, a function $f$, a
maximum number of (outer) iterations $m$, an (outer) tolerance $\varepsilon_{out}$. \\
{\bf Output: } An approximation to the leading singular triplet \\[-3mm]
\vrule height 0pt depth 0.3pt width \textwidth \\[1mm]
\begin{tabular}{rl}
1:  & Choose $\bv_1$ with $\vectornorm{\bv_1} = 1$, and set $V = [\bv_1]$, $U = \emptyset$, $M = \emptyset$, $T = \emptyset$.\\
2:  & {\bf for}         $j = 1, \dots, m$ \\
3:  & \phantom{MM}      $\bz \approx \fA\bv_{j}$ \\
4:  & \phantom{MM}      $\bu_j, \bm_j \longleftarrow \texttt{rgs}(\bz, U)$\\
5:  & \phantom{MM}      Expand basis: $U = [U, \bu_j]$\\
6:  & \phantom{MM}      Expand matrix: $M = [M, \bm_j]$ (the old $M$ is first padded with a zero row)\\
7:  & \phantom{MM}      $\bz \approx \fA^*\bu_{j}$ \\
8:  & \phantom{MM}      $\bv_{j+1},\bt_j \longleftarrow \texttt{rgs}(\bz, V)$ \\
9:  & \phantom{MM}      Expand basis: $V = [V, \bv_{j+1}]$\\
10: & \phantom{MM}      Expand matrix: $T = [T, \bt_{j}]$ \\
11: & \phantom{MM}      $K = [\texttt{zeros(j,j), M ; T , zeros(j,j)}]$\\
12: & \phantom{MM}      $[Q,D] \longleftarrow \texttt{eig}(K)$ \\
13: & \phantom{MM}      $(\theta, \bq) \longleftarrow$ with $\theta=\max_i |D_{ii}|$
 (extract $\bx,\by$ from $\bq=[\bx;\by]$ with $\|\bx\|=1,\|\by\|=1$) \\
14: & \phantom{MM}      Convergence check: if $|T(j+1,j)\bq(j)|/\theta < \varepsilon_{out}$ then
return $(\theta, \bx, \by)$ and stop \\
15: & \phantom{MM}      If required: compute variable tolerance to be used in the next iteration\\
16: & {\bf end}\\
\end{tabular}\\
\vrule height 0pt depth 0.5pt width \textwidth
}
\end{table}

\section{Practical implementation}\label{sec: Practical implementation}
Algorithm 1 implements the inexact Lanczos bidiagonalization to
approximate $\vectornorm{f(A)}_2$ and the associated singular vectors.
The function rgs$(\bz,Z)$ double orthogonalizes the vector $\bz$ with
respect to the orthogonal columns of $Z$, and returns the
orthogonalization coefficients. The same algorithm can be used to approximate
more singular triplets.

At every iteration of the Lanczos bidiagonalization, two inner iterations
(line $4$ and line $7$) approximate the corresponding matrix-vector multiplication
$f(A)\bv$ and $f(A)^*\bu$, respectively.
The inner iteration uses one of the algorithms for approximating the
action of a matrix function to a vector, as discussed in section \ref{sec:inner}.
In theory, any such algorithm could be used; in our experiments we employed
both the standard and extended Krylov subspace methods.
In case of variable inner tolerance, the next inner tolerance
is computed at the end of every Lanczos bidiagonalization iteration.
%
Moreover, after the second outer iteration,
we require that the inner stopping criterion is such that
\[
\max\{\|\bg_k^{(1)}\| , \|\bg_k^{(2)}\| \}
\le \frac{\delta_{2m,2(k-1)}}{2m\|\br_{2(k-1)}\|}\varepsilon_{out} .
\]
Note that a relative criterion is always used, that is, in practice
the quantity to be checked is divided by the current approximation $\theta^{2(k-1)}$.
This corresponds to using $\varepsilon_{out}^{(k)} = \theta^{2(k-1)} \varepsilon_{out}$ for
some fixed value $\varepsilon_{out}$.
Since $\delta_{2m,2(k-1)}$ is not available at iteration $k$, we consider the
following approximation:
\[
\delta^{2(k-1)} := \min_{\theta_{j}\in\Lambda(\Kp_{2(k-1)})\backslash \{\theta^{2(k-1)}\}} |\theta^{2(k-1)} - \theta_{j}|.
\]
In fact, $\delta_{2m,2(k-1)}$ can be much smaller than the
computed $\delta^{2(k-1)}$. However, it will not be overrated much when $\theta^{2(k-1)}$ is converging to the corresponding eigenvalue $\theta$ of $\Kp_{2m}$, since it is related to the sensitivity of $\Kp_{2m}$ and not of the matrix $\BF$. If the $\delta^{2(k-1)}$ is very small, it constrains the inner accuracy to be very small too. This occurs when the largest eigenvalues of $\Kp_{2m}$ are clustered.

\section{Numerical experiments}\label{sec:expes}
\begin{figure}[!htbp]
\centering
\includegraphics[width=1.2in,height=1.2in]{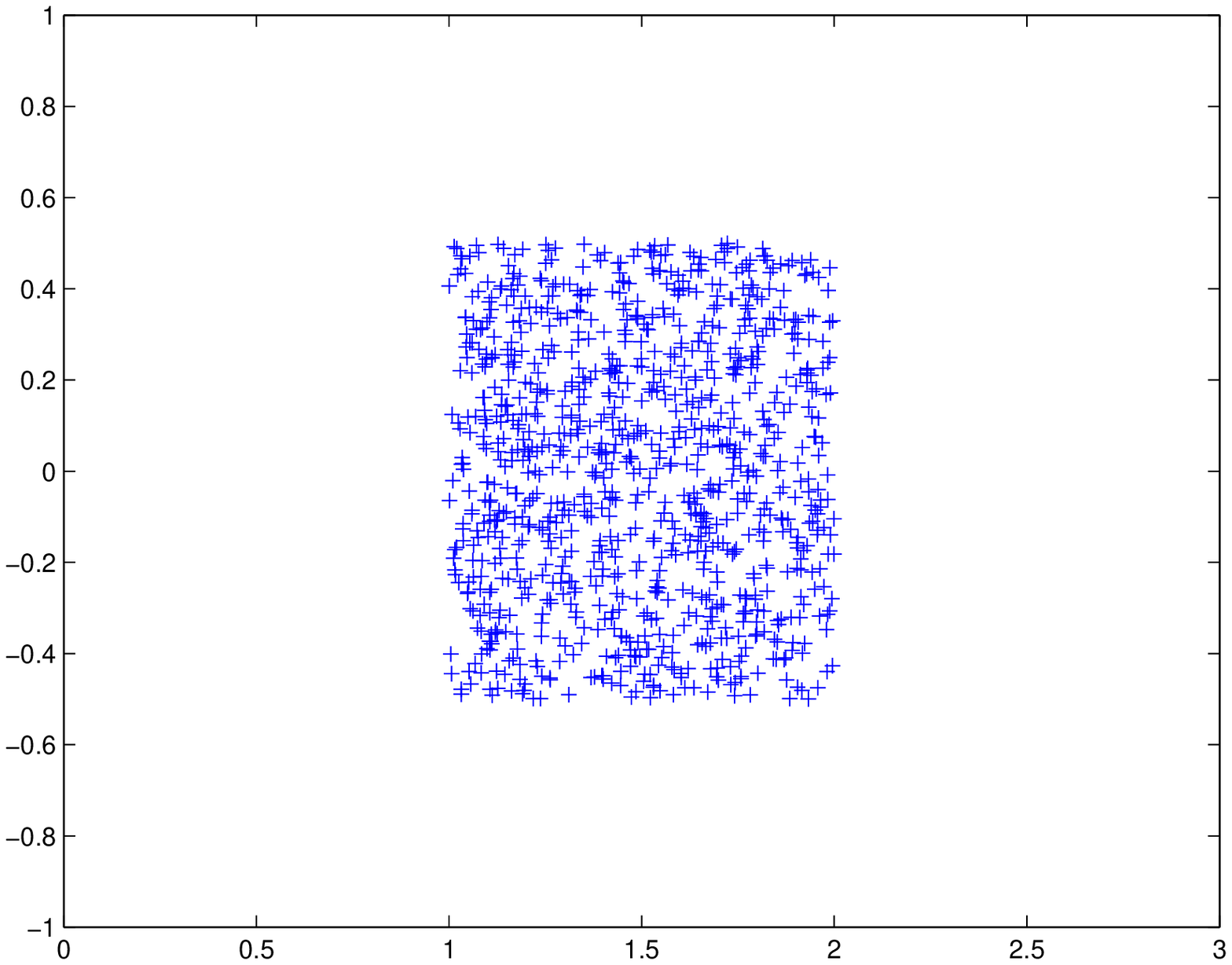}
\includegraphics[width=1.2in,height=1.2in]{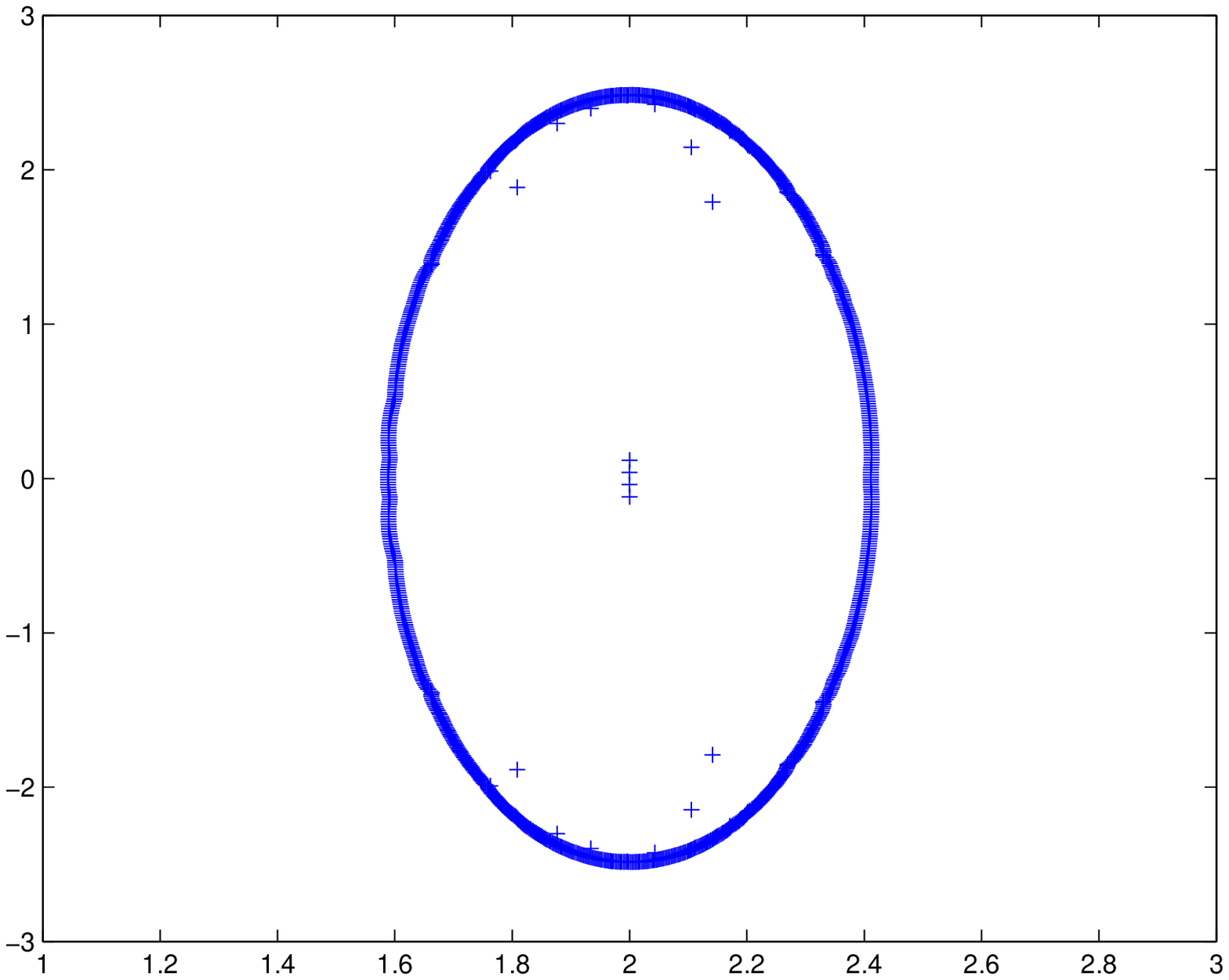}
\includegraphics[width=1.2in,height=1.2in]{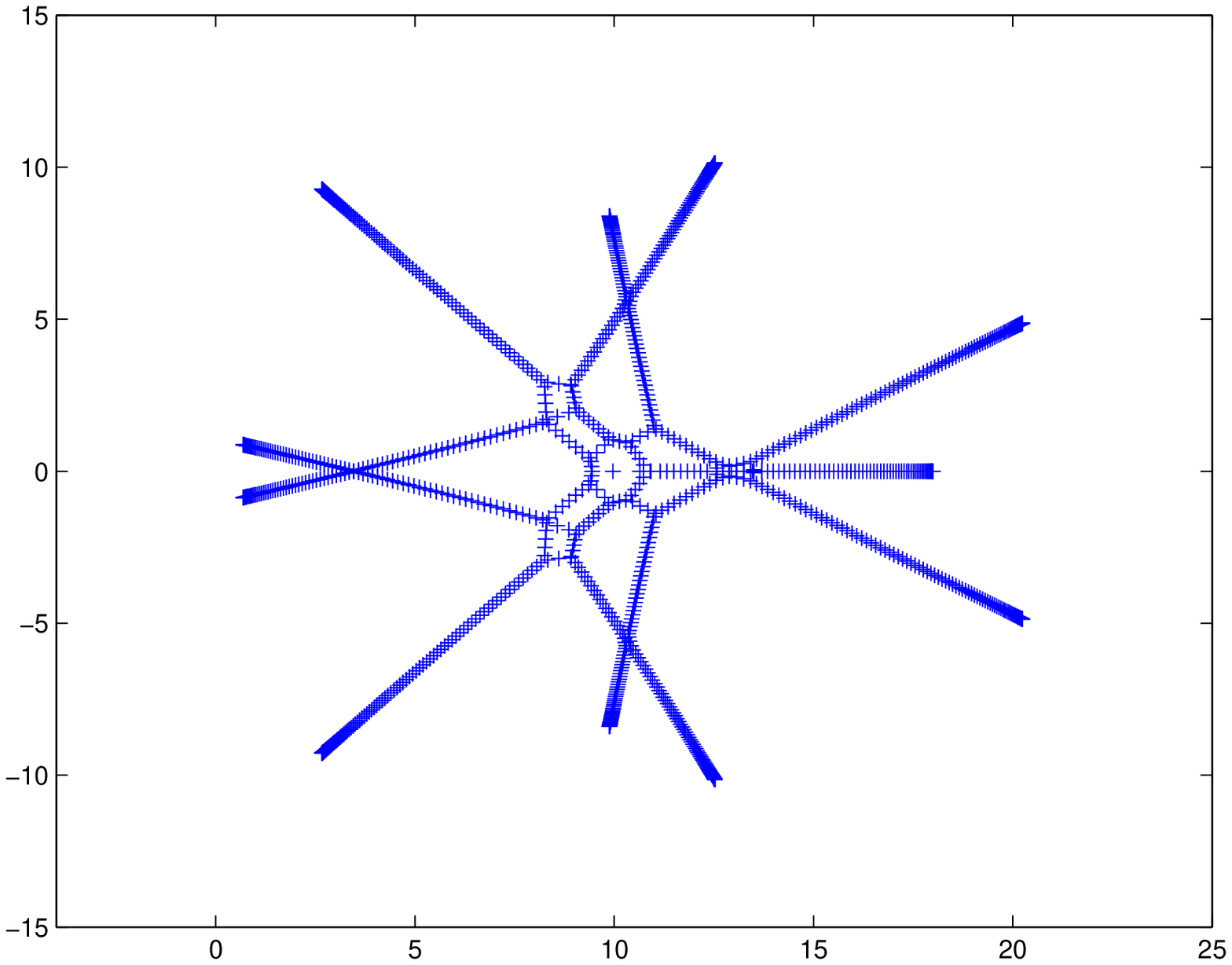}
\includegraphics[width=1.2in,height=1.2in]{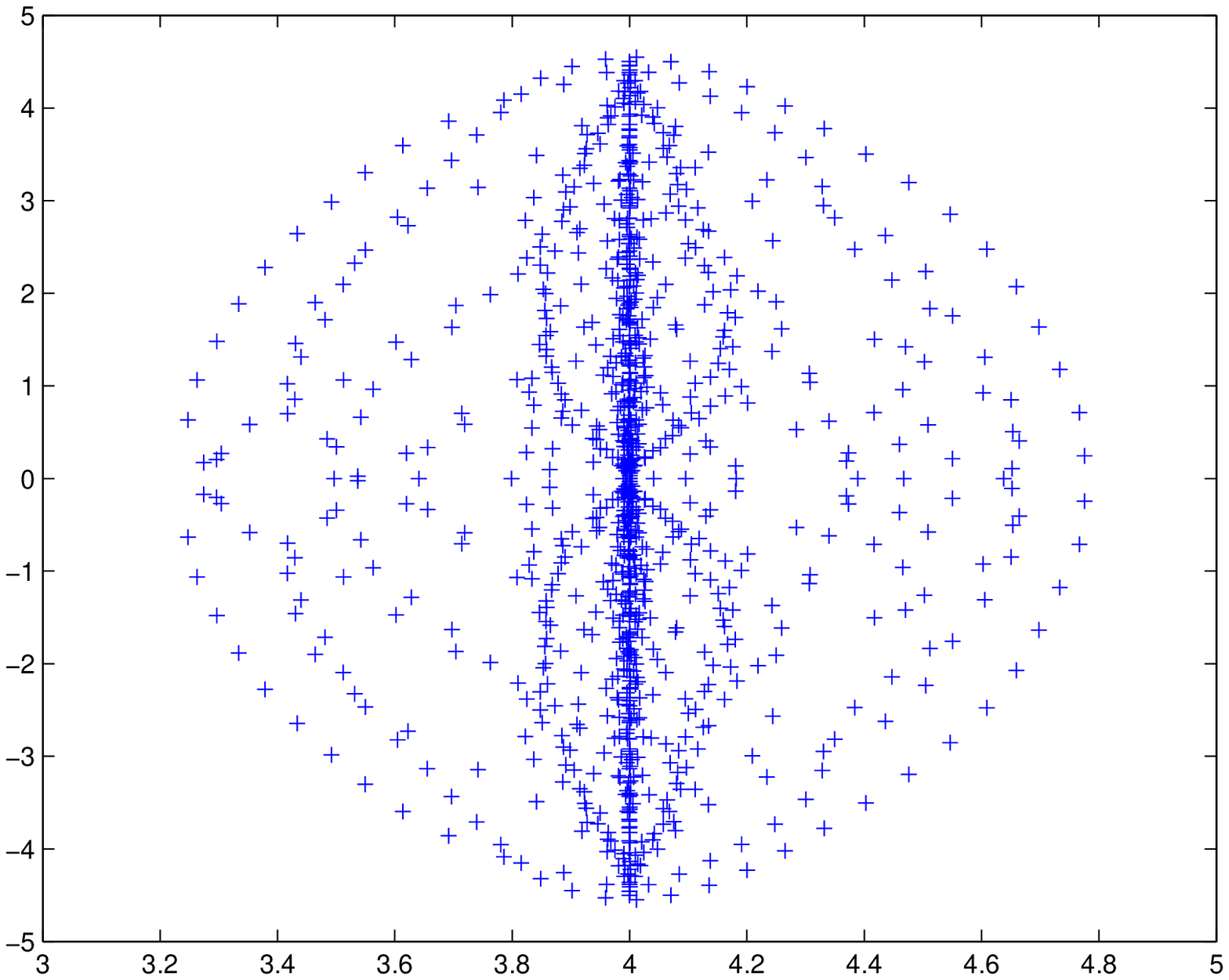}
  \caption{Spectrum of matrices $A_1, A_2, A_3, A_5$ (from left to right) in Table \ref{tab:matrices}, for
a smaller size, $n=1000$.\label{fig:spectrum}}
\end{figure}

In this section we report on our numerical experiments to evaluate the
performance of the inexact Lanczos bidiagonalization for different
combinations of matrices and functions.
All experiments were performed with Matlab Version 7.13.0.564 (R2011b) on a
Dell Latitude laptop running Ubuntu 14.04 with 4 CPUs at 2.10GHz.
We are mainly interested in the first singular triplet of $f(A)$, so as to
obtain $\|f(A)\|$.
We considered {five} different matrices, summarized in Table \ref{tab:matrices},
 all of dimension $n = 10,000$ except $A_4$. The spectrum of a sample of these matrices of
smaller size, $n=1000$, is reported in Figure~\ref{fig:spectrum}. For $A_4$, the
matrix originating from a cavity driven problem was shifted by $10 I$ so that
all functions could be treated with all methods; we refer to the Matrix Market site for
more information on this problem \cite{MatrixMarket}. For $A_5$ a
5-point stencil finite difference approximation was used, together with
homogeneous Dirichlet boundary conditions.  We considered the following functions,
$$
\exp(x), \hskip 0.2in \exp(-x), \hskip 0.2in \sqrt{x}, \hskip 0.2in \frac{1}{\sqrt{x}}, \hskip 0.2in
\frac{\exp(-\sqrt{x}) -1}{x} .
$$
We note that all these functions allow for an efficient computation
when applied to small scale matrices, by means of
specifically designed Matlab functions; see \cite{Higham2008}.
The performance also critically depends on the choice of the inner method
for approximating $f(A)\bv$ and $f(A)^*\bu$ at each iteration. We shall report
our experience with the standard and extended Krylov methods. The Rational Krylov method
could also be employed for this approximation.

\begin{table}
\centering
\begin{tabular}{|c|l|l|} \hline
Matrix & Structure & Description \\
\hline
$A_1$ & ${\rm tridiag}(0,\underline{\lambda_i}, 0.3)$ & $\lambda_i = (1+\rho_i^{(1)}) + i (\rho_i^{(2)}-0.5)$ \\
$A_2$ & ${\rm tridiag}(1.5, \underline{2}, -1)$ &  \\
$A_3$ & Toeplitz & $i$-th row:   $[4, 0, 0, 0, 0, -2, 0, \underline{10}, 0, 0, 0, 6]$ \\
$A_4$ & shifted {\sc e20r1000} & Driven cavity problem (Matrix Market) shifted as:  \\
$   $ &  &  $A :=A_{e20r1000}+10 I$ \\
$A_5$ & Sparse   & Centered Finite Difference discretization of  \\
$$ & & ${\cal L}(u) = -\nabla^2 u - 100 u_x - 100 u_y$ \\
\hline
\end{tabular}
\caption{Description of the selected matrices, all of size $n=10,000$, except $A_4$, of size 4241.
$\rho_i^{(j)}$ is a random entry taken from a uniform distribution in $(0,1)$.
\label{tab:matrices}}
\end{table}

A random vector is used to start the inexact Lanczos process.
Convergence is monitored by checking the computed residual norm with respect to
the first singular triplet, and the
inexact Lanczos bidiagonalization terminates as soon as (\ref{eqn:outer_criterion})
is satisfied; different values for $\varepsilon_{out}$ will be considered.
In case more than one triplet is desired, then all corresponding residuals
should be monitored.
In section \ref{sec:expes_bidiag} we explore the fixed inner tolerance method,
and the dependence of its performance on all the other parameters, including
the outer accuracy. Indeed, if only a rough approximation of $\|f(A)\|$ is required,
the computational efforts should proportionally low.
In section \ref{sec:expes_variable} the influence of the variable (relaxed) inner tolerance
described in section \ref{sec:variable} is
analyzed, thus  a more stringent final accuracy is considered so as to exercise
the variable inner threshold.

\begin{table}
\footnotesize
\centering
\begin{tabular}{|c| l| l| c| r| r| r| r|}\hline
matr & function & $\tilde \sigma_1$ & $\frac{\tilde\sigma_1-\tilde\sigma_2}{\tilde\sigma_1}$ & tot \#& tot \# & average & exec \\
       &   $f$  &                   &                                                        & outer  & inner & \# inner & time \\
\hline
$A_1$& $\exp(-x)$&                0.463506&3.34e-02  &   14  &  308    &   11.0    &  0.39 \\
&$\sqrt{x}$ &                     1.50143&1.29e-02  &   15  &  351    &   11.7    &  0.44 \\
&$\frac{\exp(-\sqrt{x})-1}{x}$&   0.727351&2.11e-02  &   14  &  414    &   14.8    &  0.70 \\
&$\exp(x)$&                       9.19137&1.76e-02  &   16  &  352    &   11.0    &  0.59 \\
&$1/\sqrt{x}$&                    1.10350&2.33e-02  &   13  &  364    &   14.0    &  0.73 \\
\hline
$A_2$& $\exp(-x)$&                0.222621&1.78e-02  &   11  &  308    &   14.0    &  0.24 \\
&$\sqrt{x}$ &                     1.7921&2.05e-02  &    8  &  252    &   15.8    &  0.28 \\
&$\frac{\exp(-\sqrt{x})-1}{x}$&   0.469829&1.81e-02  &   10  &  429    &   21.4    &  0.68 \\
&$\exp(x)$&                      12.1783&9.01e-02  &    5  &  139    &   13.9    &  0.14 \\
&$1/\sqrt{x}$&                    0.814356&2.26e-02  &    8  &  324    &   20.2    &  0.44 \\
\hline
$A_3$& $\exp(-x)$&                0.508086&1.48e-02  &   13  &  709    &   27.3    &  0.61 \\
&$\sqrt{x}$ &                     4.56086 &1.85e-02  &   12  & 1036    &   43.2    &  3.11 \\
&$\frac{\exp(-\sqrt{x})-1}{x}$&   0.615673&1.84e-02  &   12  & 1953    &   81.4    & 18.29 \\
&$\exp(x)$&                       6.75709$\cdot 10^8$&1.45e-02  &   13  &  694    &   26.7    &  0.83 \\
&$1/\sqrt{x}$&                    0.959018&1.70e-02  &   12  & 1852    &   77.2    & 14.56 \\
\hline
$A_4$& $\exp(-x)$&              0.000172183& 2.42e-01   &   6  &  456  &  38.0  &  0.58 \\
&$\sqrt{x}$ &                   6.09177    & 2.78e-02   &  14  &  454  &  16.2  &  0.81 \\
&$\frac{\exp(-\sqrt{x})-1}{x}$& 0.118301   & 3.07e-02   &  10  &  486  &  24.3  &  1.15 \\
&$\exp(x)$&                     3.15141$\cdot 10^{10}$ & 1.35e-01  &    5  &  397  &   39.7  & 0.49 \\
&$1/\sqrt{x}$&                  0.354039   & 5.55e-02   &   9  &  394  &  21.9  &  0.81 \\
\hline
$A_5$& $\exp(-x)$&                0.99709 & 3.39e-02 &     7 &   223  & 15.9  &   0.34 \\
&$\sqrt{x}$ &                     2.81987 & 1.16e-02 &    16 &  5145  & 160.8 & 193.35 \\
&$\frac{\exp(-\sqrt{x})-1}{x}$&   6.93384 & 2.44e-01 &     4 &  1570  & 196.2 & 111.67 \\
&$\exp(x)$&                       2959.17 & 1.84e-02 &    14 &   450  &  16.1 &   0.65 \\
&$1/\sqrt{x}$&                    7.36692 & 2.31e-01 &     4 &  1564  & 195.5 & 112.22 \\
\hline
\end{tabular}
\caption{Inexact Lanczos bidiagonalization for approximating the leading
singular triplet of $f(A)$; outer tolerance  $\varepsilon=10^{-2}$ .\label{tab:ex1e-2}}
\end{table}

\subsection{Assessing the effectiveness of the inexact bidiagonalization}\label{sec:expes_bidiag}
We analyze the performance of the inexact method when approximating
$\|f(A_i)\|$ together with the associated singular vectors. To this end, we need to monitor
the number of iterations of both the outer and the two inner iterations,
together with the execution time required to reach the required tolerance. In particular,
we show both the total and average number of inner iterations.
We also display the distance between the final first
approximate singular value, $\tilde \sigma_1$ and the second approximate
singular value, $\tilde \sigma_2$: a small relative distance implies that
the method will take more iterations to converge. Moreover, this distance
cannot be easily predicted from the matrix $A$, although it significantly influences
the computation. For instance, the largest (in modulo) eigenvalues of ${\cal F}$ associated
with the matrix function $A_2^{\frac 1 2}$ are:
\vskip 0.1in
\begin{center}
  -1.7965100, 1.7965100, 1.7964169, -1.7964169, 1.7962429, -1.7962424 .
\end{center}
\vskip 0.1in
Although this fact does not constitute a difficulty if just  the order of magnitude
of $\|A_2^{\frac 1 2}\|$ is saught, it indicates that requiring a more accurate approximation
will lead to  significantly more expensive computations.
This problem can be readily observed by comparing the outer number of iterations
in Table \ref{tab:ex1e-2} and Table \ref{tab:ex1e-4}, where we report the results of our
experiments for $\varepsilon_{out}=10^{-2}$ and
$\varepsilon_{out}=10^{-4}$, respectively. In both cases, the inner tolerance
was set to $\varepsilon_{in} = \varepsilon_{out}/(m_{\max})$, where $m_{\max}=1000$, so that
$\varepsilon_{in} =10^{-7}$ for the more stringent outer tolerance.
For all examples, the first six significant digits of $\tilde\sigma_1$ are reported.

\begin{table}
\centering
\footnotesize
\begin{tabular}{|c| l| l| c| r| r| r| r|}\hline
matr & function & $\tilde \sigma_1$ & $\frac{\tilde\sigma_1-\tilde\sigma_2}{\tilde\sigma_1}$ & tot \#& tot \# & average & exec \\
       &  $f$   &                   &                                                        & outer  & inner & \# inner & time \\
\hline
$A_1$ & $\exp(-x)$ &                0.463735&2.04e-02 &    24&    624 &      13.0 &     0.87 \\
&$\sqrt{x}$ &                       1.50496&8.76e-04 &    53&   1775 &      16.7 &     3.27 \\
&$\frac{\exp(-\sqrt{x})-1}{x}$ &    0.728200&7.62e-03 &    29&   1160 &      20.0 &     2.47 \\
&$\exp(x)$ &                        9.19576&9.22e-03 &    32&    832 &      13.0 &     1.21 \\
&$\frac{1}{\sqrt{x}}$ &             1.10504&5.52e-03 &    29&   1156 &      19.9 &     2.21 \\
 \hline
$A_2$ & $\exp(-x)$ &                0.223129&4.88e-05 &   209&   7104 &      17.0 &    26.72 \\
&$\sqrt{x}$ &                       1.79651&5.18e-05 &   162&   8069 &      24.9 &    18.92 \\
&$\frac{\exp(-\sqrt{x})-1}{x}$ &    0.470776&3.85e-05 &   193&  12320 &      31.9 &    42.03 \\
&$\exp(x)$ &                        12.1825&8.39e-04 &    47&   1596 &      17.0 &     1.41 \\
&$\frac{1}{\sqrt{x}}$ &             0.816492&5.90e-05 &   150&   9210 &      30.7 &    29.43 \\
 \hline
$A_3$ & $\exp(-x)$ &               0.509010&4.43e-05  &  224 & 14544  &    32.5  & 31.55 \\
&$\sqrt{x}$ &                      4.57175&3.35e-05  &  250 & 41844  &    83.7  & 533.17 \\
&$\frac{\exp(-\sqrt{x})-1}{x}$ &   0.616989&1.20e-04  &  155 & 39968  &   128.9 &  1078.46 \\
&$\exp(x)$ &                       6.77296$\cdot 10^8$&1.17e-04  &  183 & 11660  &   31.9  &  19.91 \\
&$\frac{1}{\sqrt{x}}$ &            0.960790&2.12e-05  &  312 & 77958  &   124.9 &  1827.25 \\
 \hline
$A_4$ & $\exp(-x)$ &               0.000172195 & 2.32e-01  &    9 &   783   &  43.5     &    1.26 \\
&$\sqrt{x}$ &                      6.09289     & 2.38e-02  &   22 &  1103   &  25.1     &    2.25 \\
&$\frac{\exp(-\sqrt{x})-1}{x}$ &   0.118347    & 2.44e-02  &   16 &  1109   &  34.7     &    3.16 \\
&$\exp(x)$ &                       3.15148$\cdot 10^{10}$ & 1.34e-01 &     7 &   626  & 44.7     &    0.87 \\
&$\frac{1}{\sqrt{x}}$ &            0.354473    & 1.18e-02  &   19 &  1227  &   32.3     &    3.91 \\
\hline
$A_5$ & $\exp(-x)$ & 0.998062&7.74e-03    & 24 &   911    & 19.0    & 1.37 \\
&$\sqrt{x}$ & 2.82811&1.67-04    &185 & 70926    & 191.7   &4059.40 \\
&$\frac{\exp(-\sqrt{x})-1}{x}$ & 6.93435 &2.32-01    &  7 &  2814    & 201.0   &186.39 \\
&$\exp(x)$ & 2975.18&2.91e-03    & 55 &  2091    & 19.0    & 3.20 \\
&$\frac{1}{\sqrt{x}}$ & 7.36768&2.17-01    &  7 &  2814    & 201.0   &197.62 \\
\hline
\end{tabular}
\caption{Inexact Lanczos bidiagonalization for approximating the leading
singular triplet of $f(A)$; outer tolerance $\varepsilon=10^{-4}$ .\label{tab:ex1e-4}}
\end{table}

Comparing the two tables also shows that the singular values are as accurate as the
outer tolerance can predict: for smaller $\varepsilon_{out}$ already the third singular
value digit changes, that is it still has to reach its final (exact) value.
This is obviously also related to the relative distance from the
second singular value, which is better captured for a smaller $\varepsilon_{out}$.

We also observe that the choice of $f$ strongly influences the overall performance:
the bidiagonalization process may take the same number of (outer) iterations for
two different selections of $f$, and yet the total computational cost be significantly
different (see $A_1$ and $A_3$ in Table \ref{tab:ex1e-2}). As a consequence,
the number of outer iterations is not a realistic measure of the algorithm complexity.

On a negative side, we observe that in both tables the method performs poorly on $A_5$ for $f(x) = \sqrt{x}$.
For this particular matrix, the inner method takes very many iterations during the whole Lanczos process, with
a number of inner iterations that is close to the average throughout. We anticipate that this is
not the case for the power method, where as the outer iterations proceed, drastically fewer
iterations are required in the inner approximation. This phenomenon seems to be
peculiar to this combination of function and matrix, since in all other cases the
performance of the Lanczos and power methods is more similar, and it will be
further investigated in a future study.

Finally, for the exponential functions $\exp(x), \exp(-x)$ we computed the upper bound
in (\ref{eqn:exp}) by using the Matlab function {\tt eigs} applied to $\frac 1 2 (A+A^*)$.
In all cases except matrix $A_4$ the estimate is pretty sharp. On the other hand,
for $\exp(A_4)$ the upper bound was $2\cdot 10^{13}$, which is three orders of magnitude
larger than the actual norm; for\footnote{This bound is obtained for $\exp(\hat A)$
with $\hat A = - A_4$.} $\exp(-A_4)$ the upper bound was 0.004737, which is more than one
order of magnitude larger than the actual value, 0.000172. This example illustrates
that, as discussed in section \ref{sec:known_methods},
the accuracy of this type of estimate cannot be easily monitored, especially in the
case of non-normal matrices.

\subsection{Comparisons with the power method}\label{sec:expes_power}
We wish to compare the performance of the new method with that of
the power method, as described in section \ref{sec:known_methods}.
Since in most cases, the leading singular values are not well isolated,
we expect that the power method will be slow if an accurate approximation
is required. Therefore, we only report results for $\varepsilon=10^{-2}$.
Moreover, our experience is that since the computation is inexact,
the product $f(A)^*(f(A)\bv)$ may give complex values, since the computed
actions of $f(A)$ and $f(A)^*$ are not the conjugate of each other. As a
result, the approximate eigenvalue  may be complex, though with a small
imaginary part, and the quantity that is actually computed is given by
$$
\lambda^{(k)} = \left | \frac{(\bv^{(k)})^* f(A)^*f(A)\bv^{(k)}}{(\bv^{(k)})^*\bv^{(k)}}\right |,
$$
where $\bv^{(k)}$ is the power method direction after $k$ iterations.
Consequently, at convergence we obtain $\tilde \sigma_1 \approx \sqrt{\lambda^{(k)}}$.
The stopping criterion is based on the relative eigenvalue residual norm, that is
$$
\|\by^{(k)} - \lambda^{(k)} \bv^{(k)}\|/\lambda^{(k)} \le \varepsilon_{out},
$$
where $\by^{(k)}$ is the result of the approximation  of $f(A)^*(f(A)\bv^{(k)})$.
Note that we kept the same tolerance as for the Lanczos bidiagonalization, although
a more stringent tolerance may be required in practice.
Table \ref{tab:power} collects the results for all test cases.

\begin{table}
\footnotesize
\centering
\begin{tabular}{|c| l| l| c| r| r| r| r|}\hline
 matr & function & tot \# &   tot \# &   $\tilde\sigma_1$ & residual &  exec \\
        &          & outer  &   inner  &                    &  norm    &  time \\
\hline
$A_1$ &$\exp(-x)$ &                 51 & 1071 & 0.46327 & 9.8648e-03 & 1.5 \\
&$\sqrt{x}$ &                       93 & 2010 & 1.5028  & 9.8607e-03 & 2.8 \\
&$\frac{\exp(-\sqrt{x})-1}{x}$ &     61 & 1782 & 0.71879 & 9.8904e-03 & 3.5 \\
&$\exp(x)$ &                        65 & 1409 & 9.1666  & 9.9964e-03 & 2.0 \\
&$1/\sqrt{x}$ &                     69 & 1942 & 1.0938  & 9.8670e-03 & 3.3 \\
\hline
$A_2$ &$\exp(-x)$ &                 36 &  899 & 0.22238 & 9.8636e-03 & 0.8 \\
&$\sqrt{x}$ &                       37 &  979 & 1.7903  & 9.7995e-03 & 1.1 \\
&$\frac{\exp(-\sqrt{x})-1}{x}$ &     36 & 1216 & 0.46921 & 9.7553e-03 & 2.1 \\
&$\exp(x)$ &                         9 &  232 & 12.176  & 9.4623e-03 & 0.2 \\
&$1/\sqrt{x}$ &                     36 & 1215 & 0.81375 & 9.8715e-03 & 1.8 \\
\hline
$A_3$ &$\exp(-x)$ &                 38 & 1605  &0.50724&  9.9024e-03&  1.5 \\
&$\sqrt{x}$ &                       41 & 1000  &4.5564 &  9.8934e-03&  1.3 \\
&$\frac{\exp(-\sqrt{x})-1}{x}$ &     34 & 4448  &0.61486&  9.9901e-03&  39.1 \\
&$\exp(x)$ &                        38 & 1699  &6.7455$\cdot 10^8$&  9.7718e-03 & 1.7 \\
&$1/\sqrt{x}$ &                     36 & 4684  &0.95774&  9.7264e-03 & 36.3 \\
\hline
$A_4$ &$\exp(-x)$ &               11&  825 &  0.00017219 & 5.8988e-03 & 1.0 \\
&$\sqrt{x}$ &                     56&  1710&  6.0870     & 9.9037e-03 & 3.0\\
&$\frac{\exp(-\sqrt{x})-1}{x}$ &   28&  1309&  0.11823    & 9.5839e-03 & 3.3\\
&$\exp(x)$ &                      10&  775 &  3.1510$\cdot 10^{10}$ & 8.8031e-03 & 1.1 \\
&$1/\sqrt{x}$ &                   33&  1361&  0.35405    & 9.9455e-03 & 2.8 \\
\hline
$A_5$ &$\exp(-x)$ &             15 &   376 &    0.99643 & 9.9168e-03 &   0.52\\
&$\sqrt{x}$ &                   52 &  1479 &    2.81329 & 9.9649e-03 &  18.47\\
&$\frac{\exp(-\sqrt{x})-1}{x}$ &  7 &  2745 &    6.93355 & 9.6566e-03 & 189.19\\
&$\exp(x)$ &                    55 &  1363 &    2956.34 & 9.8499e-03 &   1.79\\
&$1/\sqrt{x}$ &                  8 &  3137 &    7.36721 & 6.9885e-03 & 202.71\\
\hline
\end{tabular}
\caption{Power method for approximating the leading
singular triplet of $f(A)$; outer tolerance $\varepsilon=10^{-2}$\label{tab:power}}
\end{table}

As expected, the power method is more expensive than the Lanczos procedure, on
average four to five times more expensive, in
all those cases when the first singular value is not well separated from
the second one. Only for the cases of good separation, for instance with $A_5$ and the functions
$(\exp(\sqrt{x})-1)/x$ and
$1/\sqrt{x}$, convergence is reached in very few iterations, and the power
method is competitive.

We also implemented the power method as described in \cite[Algorithm 3.19]{Higham2008},
using the relative singular value residual as stopping criterion. The performance,
both in terms of inner and outer number of iterations, is comparable to that of
Table \ref{tab:power}.
Finally, we stress that in both implementations the stopping criterion
involves inexact matrix-vector products, therefore the monitored quantity
is not the true residual of the corresponding problem.

\subsection{Numerical tests with the extended Krylov subspace}\label{sec:expes_eksm}

If high accuracy is required for the final approximation, so that a more stringent
outer tolerance is used, then the inner iteration also
requires more computational effort, as its stopping tolerance is also decreased.
In this case, it may be appropriate to use more effective methods. One such possibility
is the extended Krylov subspace method \cite{Druskin.Knizhnerman.98},\cite{KnizhnermanSimoncini2010}, which may
be convenient in case the considered
function requires good approximation of the eigenvalues of $A$ closest to the origin.
In Table \ref{tab:EKSMex1e-4} we report the runs for $\varepsilon_{out}=10^{-4}$ when
EKSM is used; these numbers should be compared with those in Table \ref{tab:ex1e-4}.
We notice that EKSM requires the solution of a system with $A$ (or $A^*$) at each iteration;
to limit computational costs, a sparse LU factorization of $A$ was performed and stored once for
all at the beginning of the Lanczos bidiagonalization, and used repeatedly in the inner
iteration. This represents a tremendous saving with respect to more general rational approximations,
where solves with $(A-\tau_j I)$ have to be performed at each inner iteration, with $\tau_j$
varying with the inner step.

In Table \ref{tab:EKSMex1e-4}
all cases where EKSM provides faster convergence, that is lower execution time,
are marked in boldface. It is clear that EKSM is beneficial when good approximations to
both ends of the spectrum are required, as is the case for $x^{\alpha}$.
The lack of improvement in the case of the exponential is not unexpected, as it is known, at
least in the Hermitian case, that only one extreme of the spectrum needs to be captured
for a fast approximation of $\exp(A)\bv$.

We also remark that EKSM could also be employed as inner method in the case of the
power iteration used in section \ref{sec:expes_power}.

\begin{table}
\footnotesize
\centering
\begin{tabular}{|c| l| l| c| r| r| r| r|}\hline
matr & function & $\tilde \sigma_1$ & $\frac{\tilde\sigma_1-\tilde\sigma_2}{\tilde\sigma_1}$ & tot \#& tot \# & average & exec \\
       &          &                   &                                                        & outer  & inner & \# inner & time \\
\hline
$A_1$ & $\exp(-x)$ &              0.463735&2.04e-02   &  24   & 480    & 10.0   &   3.49  \\
&$\sqrt{x}$ &                     1.50496 &8.76e-04   &  53   & 954    & 9.0    &  8.24  \\
&$\frac{\exp(-\sqrt{x})-1}{x}$ &  0.728200&7.62e-03   &  29   & 522    & 9.0    &  4.29  \\
&$\exp(x)$ &                      9.19576 &9.22e-03   &  32   & 704    & 11.0   &   5.82  \\
&$\frac{1}{\sqrt{x}}$ &           1.10504 &5.52e-03   &  29   & 522    & 9.0    &  4.59  \\
 \hline
$A_2$ & $\exp(-x)$ &                0.223129&4.88e-05 &   209 &  5434  &    13.0 &    36.71  \\
&$\sqrt{x}$ &                       1.79651&5.18e-05  &  162  & 3564   &   11.0  &   20.03  \\
&$\frac{\exp(-\sqrt{x})-1}{x}$ &    0.470776&3.85e-05 &   193 &  4246  &     11.0&  {\bf 29.99}  \\
&$\exp(x)$ &                        12.1825&8.39e-04  &   47  & 1408   &   15.0  &    5.07  \\
&$\frac{1}{\sqrt{x}}$ &             0.816492&5.90e-05 &   150 &  3300  &     11.0&  {\bf 18.70}  \\
 \hline
$A_3$ & $\exp(-x)$ &              0.509010&4.43e-05   & 224  &11827    &   26.4  &  106.31  \\
&$\sqrt{x}$ &                     4.57175&3.35e-05    &250   &9402     & 18.8    & {\bf 84.26}  \\
&$\frac{\exp(-\sqrt{x})-1}{x}$ &  0.616989&1.20e-04   & 155  & 5578    &   18.0  &  {\bf 40.02}  \\
&$\exp(x)$ &                      6.77296$\cdot 10^8$& 1.17e-04 &   183&  11169    &   33&112.76  \\
&$\frac{1}{\sqrt{x}}$ &           0.960790&2.12e-05   & 312  &11449    &   18.3  &  {\bf 125.20}  \\
 \hline
$A_4$ & $\exp(-x)$ &              0.000172195&2.32e-01  &    9  &  376  &    20.9 &    4.20  \\
&$\sqrt{x}$ &                     6.09289   & 2.38e-02  &   22  &  483  &    11.0 &    5.99  \\
&$\frac{\exp(-\sqrt{x})-1}{x}$ &  0.118347  & 2.44e-02  &   16  &  318  &    9.9  &    4.32  \\
&$\exp(x)$ &                      3.15148$\cdot 10^{10}$ & 1.34e-01  &    7    &527   &  37.6  & 4.66  \\
&$\frac{1}{\sqrt{x}}$ &           0.354473  & 1.18e-02  &   19  &  416  &  10.9 &   4.88  \\
 \hline
$A_5$ & $\exp(-x)$ &              0.998062&7.74e-03    & 24  &  887    &   18.5  &    11.95  \\
&$\sqrt{x}$ &                     2.82811 &1.67e-04    &185  & 8165    &   22.1  &   {\bf 121.99}  \\
&$\frac{\exp(-\sqrt{x})-1}{x}$ &  6.93435 &2.32e-01    &  7  &  294    &   21.0  &    {\bf 5.01}  \\
&$\exp(x)$ &                      2975.18 &2.91e-03    & 55  & 2090    &   19.0  &   25.19  \\
&$\frac{1}{\sqrt{x}}$ &           7.36768 &2.17e-01    &  7  &  294    &   21.0  &    {\bf 4.32}  \\
\hline
\end{tabular}
\caption{Inexact Lanczos bidiagonalization, outer tolerance $\varepsilon=10^{-4}$, inner approximation: extended
Krylov subspace method.\label{tab:EKSMex1e-4}}
\end{table}

\subsection{Numerical tests with variable accuracy}\label{sec:expes_variable}
In the previous sections, for $\varepsilon_{out}=10^{-4}$
the inner tolerance was set to the fixed value
$\varepsilon_{in} =10^{-7}$. Here we explore the performance of the inexact
computation when the inner tolerance is relaxed.

A relaxed inner accuracy is most convenient when each inner iteration is
expensive, so as to profit from a lower number of inner iterations.
Therefore, we report on our experience with the extended Krylov subspace
as inner method, as the method requires one system solve with the
coefficient matrix at each iteration.
A more stringent outer tolerance was used,
that is $\varepsilon_{out} = 10^{-7}$, than in previous experiments,
so as to clearly see the relaxation in the inner tolerance; we also used
$m_{\max}=50$ as maximum number of iterations, so as to balance the
much smaller $\varepsilon_{out}$ for determining the initial
 inner tolerance.

\begin{figure}[htbp]
\centering
\includegraphics[scale=0.6]{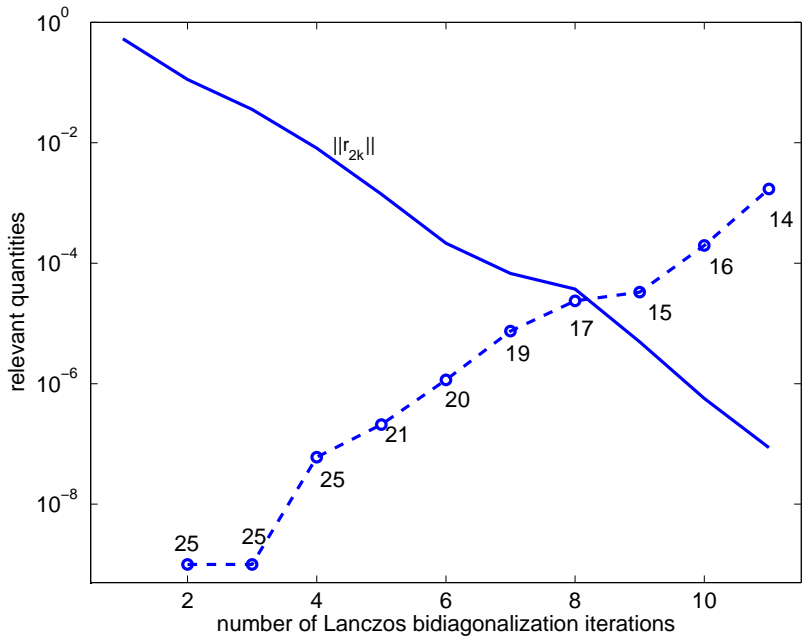}
  \caption{Relaxed inner iteration for variable stopping tolerance, to approximate
$\|A_5^{-1/2}\|$, with $\varepsilon_{out}=10^{-7}$.\label{fig:relax}}
\end{figure}

Figure \ref{fig:relax} shows the performance of the relaxation strategy
for $A_5$ and $f(x) = 1/\sqrt{x}$. The plot shows the outer convergence
history as the bidiagonalization proceeds, and the corresponding variable inner
tolerance. The digits next to each iteration report the actual numbers
of inner iterations by means of EKSM to reach the required inner accuracy
for approximating $f(A)\bv_k$; similar numbers were observed for $f(A)^*\bu_k$.

We also experimented with the approximation of more than one triplet.
We report on our findings for $A_1$ and again $f(x) = 1/\sqrt{x}$ (similar
accuracies were obtained for other functions for the same matrix); to explore
the variable inner accuracy we used $\varepsilon_{out} = 10^{-9}$ and $m_{\max}=100$.
Table \ref{tab:more_singvals} shows the largest ten singular values  obtained
with a fixed inner tolerance of $10^{-11}$ ($\tilde \sigma_j$, second column), and with
a relaxed inner tolerance ($\tilde \sigma_j^{(fl)}$, third column), which a-posteriori
we observed to go from $10^{-11}$
up to $10^{-5}$. The last column reports the relative error
$|\tilde \sigma_j - \tilde \sigma_j^{(fl)}|/\tilde \sigma_j$.
In both cases, the inexact Lanczos iteration was stopped as soon as the outer stopping criterion
was satisfied for the largest singular value. While in the fixed inner tolerance
case the number of iterations varied between 28 and 30, in the flexible case a number of
iterations as low as 15 was needed to satisfy the inner criterion at the last stage
of the convergence process. After exiting the flexible procedure, however, the first
ten approximate singular values are very close to those obtained with the fixed
inner tolerance, much closer than warranted by the final inner accuracy of $10^{-5}$.
This shows in particular that the flexible inner tolerance is conservative, and more accurate
approximations are usually expected. We refer the reader to \cite{Sim05} for a more detailed
analysis of the approximation of more than one eigenpair (and thus more than one triplet in
our context).

\begin{table}
\centering
\begin{tabular}{cccc} \hline
$j$ & $\tilde \sigma_j$ & $\tilde \sigma_j^{(fl)}$ & $|\tilde \sigma_j - \tilde \sigma_j^{(fl)}|/\tilde \sigma_j$\\
\hline
1 &1.117020718026223&1.117020718026212&9.939144428645173e-15 \\
2 &1.107884805324699&1.107884805324724&2.264769787972212e-14 \\
3 &1.098394607515649&1.098394607513931&1.564063676705998e-12 \\
4 &1.095557563655289&1.095557550135225&1.234080654983002e-08 \\
5 &1.087939266226247&1.087939266157844&6.287397018996839e-11 \\
6 &1.081739455193175&1.081739454786304&3.761265981356762e-10 \\
7 &1.077326677541678&1.077326677826174&2.640758724893079e-10 \\
8 &1.070641401649297&1.070641400153385&1.397211109700097e-09 \\
9 &1.064637797334345&1.064637795718615&1.517633507177347e-09 \\
10&1.055679471834666&1.055679470916211&8.700132135409333e-10 \\
\hline
\end{tabular}
\caption{First ten approximate singular values of $A_1^{-1/2}$ with fixed tolerance
($\varepsilon_{out}=10^{-9}$), and relaxed inner tolerance.\label{tab:more_singvals}}
\end{table}


\section{Final considerations}\label{sec: discussion}
We have explored the use of an inexact Lanczos bidiagonalization method
for approximating the leading singular triplet of a large matrix function, and
in particular its 2-norm. Although several strategies are known to
provide rough estimates of a matrix function 2-norm, more accurate
approximations require a careful implementation of available approaches,
since neither $f(A)$ nor products of the type $f(A)\bv$ are available exactly.
In particular, we showed that the Lanczos bidiagonalization yields a non-Hermitian
perturbation of the original Hermitian matrix, and the recurrence needs
to be revisited.
%
Our numerical experiments showed that the computational complexity may vary
significantly depending on the requested final accuracy, since
the two inner iterations for the approximation of $f(A)\bv$ and
$f(A)^*\bu$ may be very time and memory consuming.
We showed that the relaxed strategy alleviates this problem whenever accurate approximations
are required. However, for particular selections of matrices and functions,
the approximation of $f(A)\bv$ can still be very expensive, and some
other strategies could be exploited, such as restarting; see, e.g.,
 \cite{Eiermann2006},\cite{FrommerGuettelSchweitzer.14},\cite{Guettel.survey.13} and references therein.
Finally, our approach could be used to estimate the norm of other
matrix objects, such as the geometric mean \cite{Bhatia.Grover.12},
or the {\it derivatives} of matrix functions, such as
the Fr{\'e}chet derivative of the matrix exponential or of other functions \cite{Higham.Relton.14}.


\begin{thebibliography}{10}

\bibitem{Baxter.94}
{\sc B.~J.~C. Baxter}, {\em Norm estimates for inverses of {T}oeplitz distance
  matrices}, J. Approximation theory, 79 (1994), pp.~222--242.

\bibitem{BenziSIREV.13}
{\sc M.~Benzi, P.~Boito, and N.~Razouk}, {\em Decay properties of spectral
  projectors with applications to electronic structure}, SIAM Review, 55
  (2013), pp.~3--64.

\bibitem{Bhatia.Grover.12}
{\sc R.~Bhatia and P.~Grover}, {\em Norm inequalities related to the matrix
  geometric mean}, Lin. Alg. Appl., 437 (2012), pp.~726--733.

\bibitem{BoF00}
{\sc A.~Bouras and V.~Frayss\'e}, {\em Inexact matrix-vector products in {
  Krylov} methods for solving linear systems: A relaxation strategy}, SIAM J.
  Matrix Analysis and Appl., 26 (2005), pp.~660 -- 678.

\bibitem{Choi.PhD.2013}
{\sc D.~Choi}, {\em Estimating Norms of Matrix Functions using Numerical
  Ranges}, PhD thesis, University of Washington, 2013.

\bibitem{Crouzeix2007}
{\sc M.~Crouzeix}, {\em {Numerical range and functional calculus in Hilbert
  space}}, J.\ Functional Analysis, 244 (2007), pp.~668--690.

\bibitem{Cul02}
{\sc J.~K. Cullum and R.~A. Willoughby}, {\em Lanczos Algorithms for Large
  Symmetric Eigenvalue Computations: Vol. 1: Theory}, vol.~41, SIAM, 2002.

\bibitem{Druskin.Knizhnerman.98}
{\sc V.~Druskin and L.~Knizhnerman}, {\em Extended {Krylov} subspaces:
  approximation of the matrix square root and related functions}, SIAM J.
  Matrix Anal. Appl., 19 (1998), pp.~755--771.

\bibitem{druskin:3760}
{\sc V.~Druskin, L.~Knizhnerman, and M.~Zaslavsky}, {\em Solution of large
  scale evolutionary problems using rational krylov subspaces with optimized
  shifts}, Tech. Rep.~5, 2009.

\bibitem{Eiermann2006}
{\sc M.~Eiermann and O.~Ernst}, {\em A restarted {Krylov} subspace method for
  the evaluation of matrix functions}, SIAM J. Numer. Anal., 44 (2006),
  pp.~2481--2504.

\bibitem{FrommerGuettelSchweitzer.14}
{\sc A.~Frommer, S.~G{\"u}ttel, and M.~Schweitzer}, {\em Convergence of
  restarted {K}rylov subspace methods for {S}tieltjes functions of matrices},
  SIAM J. Matrix Anal. Appl,  (2014).

\bibitem{FrommerMarch2006}
{\sc A.~Frommer and V.~Simoncini}, {\em Matrix functions}, tech. rep.,
  Dipartimento di Matematica, Bologna, I, March 2006.
\newblock To appear on SIAM J. Scient. Computing.

\bibitem{Gear.81}
{\sc C.~W. Gear}, {\em Numerical solution of ordinary differential equations:
  Is there anything left to do?}, SIAM Rev., 23 (1981), pp.~10--24.

\bibitem{Gil.2010}
{\sc M.~Gil'}, {\em Perturbation of functions of diagonalizable matrices},
  Electronic J. of Linear Algebra, 20 (2010), pp.~303--313.

\bibitem{GLo96}
{\sc G.~H. Golub and C.~F. {Van Loan}}, {\em Matrix {C}omputations}, The John
  Hopkins University Press, Baltimore, London, 3rd~ed., 1996.

\bibitem{Gustafson1997}
{\sc K.~Gustafson and D.~K.~M. Rao}, {\em Numerical Range: The Field of Values
  of Linear Operators and Matrices}, Universitext, {Springer}, New York, NY,
  USA, 1997.

\bibitem{Guettel.survey.13}
{\sc S.~G{\"u}ttel}, {\em Rational {K}rylov approximation of matrix functions:
  Numerical methods and optimal pole selection}, GAMM Mitteilungen, 36 (2013),
  pp.~8--31.

\bibitem{Higham.Relton.14}
{\sc N.~Higham and S.~Relton}, {\em Estimating the condition number of the
  frechet derivative of a matrix function}, SIAM J. Sci. Comput.,  (2014).

\bibitem{Higham2008}
{\sc N.~J. Higham}, {\em Matrix Functions -- Theory and Applications}, SIAM,
  Philadelphia, USA, 2008.

\bibitem{Hochbruck1999}
{\sc M.~Hochbruck and C.~Lubich}, {\em {Exponential integrators for
  quantum-classical molecular dynamics}}, BIT, Numerical Mathematics, 39
  (1999), pp.~620--645.

\bibitem{Hochbruck.Ostermann.10}
{\sc M.~Hochbruck and A.~Ostermann}, {\em Exponential integrators}, Acta
  Numerica, 19 (2010), pp.~209--286.

\bibitem{Horn.Johnson.91}
{\sc R.~A. Horn and C.~R. Johnson}, {\em {Topics in Matrix Analysis}},
  Cambridge University Press, Cambridge, 1991.

\bibitem{Horn.Johnson.13}
\leavevmode\vrule height 2pt depth -1.6pt width 23pt, {\em {Matrix Analysis}},
  Cambridge University Press, Cambridge, {II}~ed., 2013.

\bibitem{Jarlebring.Rubensson.12}
{\sc E.~Jarlebring and E.~H. Rubensson}, {\em On the condition number and
  perturbation of matrix functions for hermitian matrices}.
\newblock arXiv:1206.1762v1, June 2012.

\bibitem{KnizhnermanSimoncini2010}
{\sc L.~Knizhnerman and V.~Simoncini}, {\em {A new investigation of the
  extended Krylov subspace method for matrix function evaluations}}, Numerical
  Linear Algebra with Applications, 17 (2010), pp.~615--638.

\bibitem{Larsen.tr98}
{\sc R.~M. Larsen}, {\em Lanczos bidiagonalization with partial
  reorthogonalization}, Tech. Rep. DAIMI PB-357, Department of Computer
  Science, Aarhus University, 1998.

\bibitem{MatrixMarket}
{\sc {Matrix Market}}, {\em A visual repository of test data for use in
  comparative studies of algorithms for numerical linear algebra,
  {M}athematical and {C}omputational {S}ciences {D}ivision, {N}ational
  {I}nstitute of {S}tandards and {T}echnology}.
\newblock Online at {\tt http://math.nist.gov/MatrixMarket}.

\bibitem{Mengi.Overton.05}
{\sc E.~Mengi and M.~Overton}, {\em Algorithms for the computation of the
  pseudospectral radius and the numerical radius of a matrix}, IMA J. Num. An.,
  25 (2005), pp.~648--669.

\bibitem{Rubensson.12}
{\sc E.~H. Rubensson}, {\em Controlling errors in recursive {Fermi-Dirac}
  operator expansions with applications in electronic structure theory}, SIAM
  J. Sci. Comput., 34 (2012), pp.~B1--B23.

\bibitem{Schmitt.83}
{\sc B.~Schmitt}, {\em Norm bounds for rational matrix functions}, Numerische
  Mathematik, 42 (1983), pp.~379--389.

\bibitem{Sim05}
{\sc V.~Simoncini}, {\em Variable accuracy of matrix-vector products in
  projection methods for eigencomputation}, SIAM J. Numer. Anal., 43 (2005),
  pp.~1155--1174.

\bibitem{Simoncini2003c}
{\sc V.~Simoncini and D.~B. Szyld}, {\em Theory of inexact {Krylov} subspace
  methods and applications to scientific computing}, SIAM J. Sci. Comput., 25
  (2003), pp.~454--477.

\bibitem{StS90}
{\sc G.~W. Stewart and J.-G. Sun}, {\em Matrix perturbation theory}, Computer
  Science and Scientific Computing, Academic Press Inc., Boston, MA,  (1990).

\bibitem{Trefethen.Embree.05}
{\sc L.~Trefethen and M.~Embree}, {\em Spectra and pseudospectra. The behavior
  of non-normal matrices and operators}, Princeton University Press, 2005.

\bibitem{Watson.96}
{\sc G.~A. Watson}, {\em Computing the numerical radius}, Lin. Alg. Appl., 234
  (1996), pp.~163--172.

\end{thebibliography}

\appendix
\section{Proof of Proposition \ref{prop:tau}}
Define the submatrix of $\Kp_{2m}$ of size $2k$ as
$$
\Kp_{2k} = \begin{bmatrix}0 & M_k \\ T_k & 0 \end{bmatrix} ,
\qquad {\rm i.e.} \qquad
\Kp_{2m} = \begin{bmatrix}
0   & 0        & M_k  & M_{\star}\\
0   & 0        & 0    & \star \\
T_k                   & T_{\star} & 0    & 0\\
t_{k+1,k}\be_1\be_k^* & \star     & 0    & 0\\
\end{bmatrix} .
$$
Define the vector $\tilde{\bq} = \frac{1}{\sqrt{2}}[\bx; 0; \by; 0]$, where the $0$-vectors have length $m-k$.
Let $\mathcal{X} = [\tilde{\bq}, Y]$
be such that $\mathcal{X}$ is unitary,
where $Y = [Y_1; Y_2; Y_3; Y_4]$.
This implies that $\frac{1}{\sqrt{2}}(Y_1^*\bx + Y_3^*\by) = 0$, $Y_4Y_4^* = I = Y_2Y_2^*$ and $Y_2Y_4^* = Y_4Y_2^* = 0$. Now, write
\[
\mathcal{X}^*\Kp_{2m}\mathcal{X} =
\begin{bmatrix}
\tilde{\bq}^*\Kp_{2m}\tilde{\bq} & \tilde{\bq}^*\Kp_{2m}Y\\
Y^*\Kp_{2m}\tilde{\bq}   & Y^*\Kp_{2m}Y\\
\end{bmatrix}
=
\begin{bmatrix}
\theta^{(2k)} & \bg_1^*\\
\bg_2  & \underline{\Kp}_{2m}\\
\end{bmatrix} .
\]
Here
\[\|\bg_2\| = \|Y^*\Kp_{2m}\tilde{\bq}\| =
\| \frac{1}{\sqrt{2}}Y^*\begin{bmatrix}M_k\by\\ 0 \\ T_k\bx \\ t_{k+1,k}\be_1\be_k^*\bx\end{bmatrix}\|
= \|\frac{1}{\sqrt{2}}Y_4^*t_{k+1,k}\be_1\be_k^*\bx\| = \|\br_{2k}\|.
\]
Further, since $\bs_{2m}^*Y = \tilde{\bq}^*\Kp_{2m}Y - \theta^{(2k)}\tilde{\bq}^*Y = \tilde{\bq}^*\Kp_{2m}Y$, we have
\[\|\bg_1\| = \|\tilde{\bq}^*\Kp_{2m}Y\| = \|\bs_{2m}^*Y\| \le \|\bs_{2m}\|.
\]
Now, by \cite[Theorem 2.1, p.230]{StS90},
\[
\text{if}\,\,\,\frac{\|\br_{2k}\|\,\|\bs_{2m}\|}{\delta_{2m,2k}^2} < \frac{1}{4}, \quad \text{i.e.,}\quad \|\br_{2k}\| < \frac{\delta_{2m,2k}^2}{4\|\bs_{2m}\|},
\]
then there exists a vector $\bp\in\mathcal{C}^{2m-1}$ satisfying $\tau = \|\bp\| < 2 \frac{\|\br_{2k}\|}{\delta_{2m,2k}}$, such that the unit norm vector
$$
\bq =
\begin{bmatrix}\bx_1 \\ \bx_2 \\ \by_1 \\ \by_2 \end{bmatrix} =
\frac{1}{\sqrt{1+\|\bp\|^2}}
\left(
\frac{1}{\sqrt{2}}\begin{bmatrix}\bx \\ 0 \\ \by \\ 0 \end{bmatrix}
+ \begin{bmatrix}Y_1 \\ Y_2 \\ Y_3 \\ Y_4 \end{bmatrix} \bp \right )
$$
is an eigenvector of $\Kp_{2m}$. Moreover,
\[
\vectornorm{\begin{bmatrix}\bx_2 \\ \by_2 \end{bmatrix}}
= \vectornorm{\frac{1}{\sqrt{1+\tau^2}}\begin{bmatrix}Y_2 \\ Y_4 \end{bmatrix}\bp}
\le \frac{\tau}{\sqrt{1+\tau^2}} .
\]
Further, this same theorem states that
$|\theta - \theta^{(2k)}| = \|\bg_1^*\bp\| \le \|\bg_1\| \|\bp\| \le \|\bs_{2m}\|\tau$. \endproof


\end{document}